\documentclass[review,hidelinks,onefignum,onetabnum]{siamart220329}
\usepackage[T1]{fontenc}
\usepackage{amsfonts,amssymb,mathtools,bm}
\usepackage{graphicx}
\usepackage{booktabs}
\usepackage{placeins}
\usepackage{enumitem}
\usepackage{stmaryrd}
\usepackage{xcolor}
\newsiamremark{remark}{Remark}
\newsiamthm{assumption}{Assumption}
\newsiamthm{example}{Example}

\makeatletter
\def\refstepcounter@optarg[#1]#2{\cref@old@refstepcounter{#2}\cref@constructprefix{#2}{\cref@result}\@ifundefined{cref@#1@alias}{\def\@tempa{#1}}{\def\@tempa{\csname cref@#1@alias\endcsname}}\protected@edef\cref@currentlabel{[\@tempa][\arabic{#2}][\cref@result]\csname p@#2\endcsname\csname the#2\endcsname}}
\makeatother

\newcommand{\paramtablerules}{\color{black}\setlength{\heavyrulewidth}{0.8pt}\setlength{\lightrulewidth}{0.5pt}\setlength{\cmidrulewidth}{0.3pt}\setlength{\arrayrulewidth}{0.5pt}}

\newcommand{\papertablestyle}{\scriptsize\renewcommand{\arraystretch}{0.92}\setlength{\tabcolsep}{1.5pt}}

\newcommand{\bs}{\bm}    \newcommand{\Ac}{\mathcal A} \newcommand{\tI}{\underline I}  \newcommand{\tauh}{\underline\tau_h}  \newcommand{\uu}{\bm u_h}      \newcommand{\gh}{\gamma_h}  \newcommand{\normal}{\bm n} \newcommand{\jump}[1]{\left\llbracket #1\right\rrbracket}  \newcommand{\norm}[1]{\left\lVert #1\right\rVert} \newcommand{\abs}[1]{\left\lvert #1\right\rvert} \newcommand{\tr}{\operatorname{tr}} \newcommand{\asop}{\operatorname{as}}  \newcommand{\diag}{\operatorname{diag}} \newcommand{\Div}{\operatorname{div}}  \newcommand{\tcurl}{\operatorname{curl}} \newcommand{\trot}{\operatorname{rot}} 

       \newcommand{\Mp}{\mathsf M_p} \newcommand{\App}{\mathsf A_{pp}}

\newcommand{\A}{\mathsf A} \newcommand{\B}{\mathsf B} \newcommand{\C}{\mathsf C} \newcommand{\D}{\mathsf D} \newcommand{\E}{\mathsf E} \newcommand{\F}{\mathsf F} \newcommand{\M}{\mathsf M}

\newcommand{\Kuu}{\mathsf K_{uu}} \newcommand{\Kug}{\mathsf K_{u\gamma}} \newcommand{\Kup}{\mathsf K_{up}} \newcommand{\Kggzero}{\mathsf K_{\gamma\gamma}^{0}} \newcommand{\Jg}{\mathsf J_{\gamma}} \newcommand{\Mg}{\mathsf M_{\gamma}} \newcommand{\Kgg}{\mathsf K_{\gamma\gamma}} \newcommand{\Kgp}{\mathsf K_{\gamma p}} \newcommand{\Kmzero}{\mathsf K_m^{0}} \newcommand{\Km}{\mathsf K_m} \newcommand{\Dm}{\mathsf D_m} \newcommand{\Lop}{\mathsf L} \newcommand{\Spsig}{\mathsf S_p^\sigma} \newcommand{\Spz}{\mathsf S_p^z} \newcommand{\Tp}{\mathsf T_p}  \newcommand{\Pmz}{\mathsf P_m^z} \newcommand{\Sp}{\mathsf S_{p}} \newcommand{\DT}{\mathcal D_{\Tp}}  \newcommand{\DPmz}{\mathcal D_{\Pmz}}

\newcommand{\ssig}{\boldsymbol\sigma} \newcommand{\uc}{\mathbf u} \newcommand{\gammac}{\boldsymbol\gamma} \newcommand{\zc}{\widehat{\mathbf z}} \newcommand{\pp}{\mathbf p}\newcommand{\qq}{\mathbf q} \newcommand{\mm}{\mathbf m} \newcommand{\nn}{\boldsymbol{\zeta}}

\newcommand{\Ared}{\mathcal A_{\rm red}}  \newcommand{\Psc}{\mathsf P_{\rm sc}}        \newcommand{\spec}{\operatorname{spec}}

\newcommand{\Mzero}{\mathsf M_0} \newcommand{\Kdr}{K_{\rm dr}} \newcommand{\clambda}{c_{\lambda}} \newcommand{\Thetah}{\Theta_h} \newcommand{\Gammah}{\mathfrak G_h}

\title{Parameter-robust preconditioners for a stabilized control-volume discretization of Biot's equations}
\author{
	Lei Wei\thanks{Department of Mathematics, City University of Hong Kong, Kowloon Tong, Hong Kong SAR, China. Email: \email{leiwei@cityu.edu.hk}.}
	\and
	Lina Zhao\thanks{Corresponding author. Department of Mathematics, City University of Hong Kong, Kowloon Tong, Hong Kong SAR, China. Email: \email{linazha@cityu.edu.hk}. {L. Zhao is partially supported by Research Grants Council of the Hong Kong Special Administrative Region, China (Project No. CityU  11302424, CityU 11303225)
	and Young Scientists fund, NSFC (Project No. 12401543).}
}
}
\headers{Parameter-robust preconditioners for Biot}{L. Wei and L. Zhao}

\begin{document}
	\nolinenumbers
	\maketitle

	\begin{abstract}
	We develop and analyze parameter-robust block preconditioners for a stabilized control-volume discretization of the quasi-static Biot equations. Beginning with a five-field formulation, we eliminate the poroelastic stress and time-scaled Darcy velocity through independent vertex-local solves, obtaining a symmetric cell-centered displacement--rotation--pressure system. A macrocell rotation-jump term stabilizes the weak-symmetry constraint. We prove the well-posedness of the fully discrete scheme and show that the condensed, coupled displacement--rotation mechanics block is uniformly spectrally equivalent to its block-diagonal counterpart. For uncoupled preconditioners, we derive an exact discrete pressure-control criterion and a sufficient parameter condition that accommodates the scaling \(c_0\asymp K_{\rm dr}^{-1}\), thereby allowing the storage coefficient to vanish in the nearly incompressible limit. We then construct a Schur-consistent pressure block that approximates the mechanics-induced pressure correction and prove its uniform equivalence to the exact pressure Schur complement.  These results yield parameter-independent spectral bounds for block-diagonal preconditioning and weighted field-of-values estimates for lower-triangular preconditioning. Practical realizations employ a macrocell mass-matrix approximation for the stabilized rotation block and tailored two-grid solvers for the displacement and pressure blocks. Numerical experiments verify the expected convergence behavior in two dimensions and demonstrate robust performance under mesh refinement and broad variations of the physical parameters in both two and three dimensions.
	\end{abstract}



	\section{Introduction}	\label{sec:introduction}

	The quasi-static Biot model couples porous-solid deformation and pore-fluid flow, with applications in geomechanics and subsurface flow \cite{Biot1941,Coussy2004}. Robust discretizations and solvers must handle nearly incompressible solids, small or vanishing storage, and weak Darcy control caused by low permeability or small time steps. These challenges motivate parameter-robust formulations and block preconditioners based on parameter-weighted norms \cite{MardalWinther2011,LeeMardalWinther2017}.

	For Biot's equations, existing discretizations include classical displacement--pressure formulations \cite{MuradLoula1994,AdlerEtAl2018}, locking-free total-pressure formulations \cite{OyarzuaRuizBaier2016,LeeMardalWinther2017}, and fully mixed formulations retaining stress and Darcy velocity \cite{Lee2016,AmbartsumyanKhattatovYotov2020}. The latter directly approximate both fields and, with suitable $H(\operatorname{div})$-conforming spaces, preserve local momentum and mass balances. For plane elasticity, Arnold et al.~\cite{ArnoldDouglasGupta1984}
	established mixed finite element error estimates uniform with respect to
	compressibility, including the incompressible limit.
	Weak stress symmetry is imposed through a rotation multiplier \cite{ArnoldFalkWinther2007}. Before condensation, fully mixed formulations lead to large four- or five-field saddle-point systems. Local condensation and hybridization reduce the globally coupled unknowns while retaining conservation properties \cite{Nordbotten2016,KrausEtAl2021}.

    Multipoint flux and stress methods solve local systems on interaction regions formed by cells or subcells sharing a mesh vertex \cite{Aavatsmark2002,WheelerYotov2006,AmbartsumyanEtAlSimplicial2019,AmbartsumyanEtAlQuadrilateral2019}. In these multipoint mixed finite element methods, vertex quadrature produces block-diagonal flux or stress matrices, enabling local elimination around each vertex.  The closely related coupled multipoint stress--multipoint flux mixed finite element method for Biot's equations \cite{AmbartsumyanKhattatovYotov2020} uses lowest-order BDM stress and Darcy velocity, cellwise constant displacement and pressure, and continuous piecewise linear or bilinear rotations. Quadrature of the velocity, stress, and stress--rotation forms permits local elimination of stress and Darcy velocity followed by rotation, yielding a cell-centered displacement--pressure matrix that is block-skew-symmetric and positive definite.
	Fu and Zhao \cite{FuZhao2025CiCP} instead developed a locking-free control-volume method for elasticity on quadrilateral grids, with local momentum conservation and piecewise constant stress, displacement, and rotation that simplify implementation. Their stress space localizes the compliance form to vertex interaction regions by construction, eliminating the need for vertex quadrature; in the macrocell-rotation variant, the globally coupled unknowns are the cell-centered displacement and rotation. We extend this variant to Biot's equations using the same subcell-wise constant normal-flux space for stress rows and Darcy velocity, and macrocell-wise constant displacement, rotation, and pressure. The stress and scaled Darcy blocks are block diagonal over vertex interaction regions by construction and are eliminated independently, yielding a symmetric displacement--rotation--pressure system. A macrocell rotation-jump stabilization acts only in the weak-symmetry equation, preserving local momentum and mass balances.

Parameter-robust block preconditioners for Biot systems have been developed primarily within the operator-preconditioning framework and through closely related norm-based approaches. Operator preconditioning equips the underlying product space with parameter-weighted norms for which the variational operator satisfies parameter-independent continuity and inf--sup estimates. The inverse Riesz map associated with these norms, or a uniformly spectrally equivalent approximation, then defines a robust preconditioner \cite{Hiptmair2006,MardalWinther2011}. This framework has been applied to total-pressure formulations \cite{LeeMardalWinther2017,Lee2023} and to conservative four-field total-pressure formulations \cite{BoonEtAl2021}. Causemann and Kuchta~\cite{CausemannKuchta2026SISC} developed a scalable, parameter-robust solver for cell-by-cell poroelasticity based on norm-equivalent preconditioning and algebraic multigrid approximations of the diagonal blocks. Related norm-based constructions have been proposed for classical two- and three-field formulations \cite{ChenEtAl2020}, stabilized discretizations \cite{AdlerEtAl2020}, and hybridized systems \cite{KrausEtAl2021,HenriquezEtAl2025Biot}. Huang and Wang~\cite{HuangWang2026SISC} introduced parameter-robust block-triangular preconditioners for two- and three-field poroelasticity formulations that avoid the explicit construction of a Schur complement.
For statically condensed systems, robustness must be established for the reduced operator itself: one may reduce a full-system preconditioner under suitable uniform conditions \cite{HenriquezEtAl2025SISC,HenriquezEtAl2025Biot} or approximate the condensed blocks directly \cite{FuKuang2023}. We adopt the latter approach for the reduced displacement--rotation--pressure system. Exact Schur factorizations provide the algebraic benchmark for the mechanics-induced pressure correction, while weighted field-of-values (FOV) estimates are used to analyze the nonsymmetric triangularly preconditioned operators \cite{MurphyGolubWathen2000,Ipsen2001,LoghinWathen2004}.

	Two issues must be resolved to construct robust preconditioners for the condensed system. First, separating displacement and rotation makes mechanics solves cheaper, but the compliance operator loses uniform coercivity on spherical stresses near incompressibility; uniform control of the displacement--rotation coupling is therefore needed. Second, when storage and Darcy control are weak, the natural pressure block may fail to control the mechanics--pressure coupling. The pressure-control criterion for an uncoupled strategy must therefore account jointly for storage, residual compliance, and Darcy contributions, with a bound independent of the mesh size and physical parameters. At zero storage, this strategy remains valid only if the natural pressure block is positive definite and the criterion holds uniformly. A Schur-consistent alternative approximates the mechanics-induced pressure correction explicitly. Our analysis distinguishes these strategies and their applicable parameter regimes. The main contributions are summarized as follows.
	\begin{enumerate}[label=(\roman*),leftmargin=2.4em]
		\item
		We establish the well-posedness of the proposed fully discrete control-volume scheme for the Biot model using inf--sup estimates for stabilized linear elasticity and Darcy flow.
		Notably, the rotation-jump stabilization added to the weak-symmetry equation enables the elasticity inf--sup estimate and yields a condensed rotation block uniformly spectrally equivalent to a \(\mu\)-scaled mass matrix.

		\item We prove spectral equivalence between the coupled displacement--rotation block and its block-diagonal part.
		For uncoupled preconditioners, we establish an exact discrete pressure-control criterion admitting \(c_0\asymp\Kdr^{-1}\), so storage may vanish as \(\lambda\to\infty\).

		\item We construct a computable Schur-consistent pressure block and prove its uniform spectral equivalence to the exact pressure Schur complement, yielding parameter-independent spectral bounds for block-diagonal preconditioning. We also establish uniform weighted FOV bounds for the exact and Schur-consistent lower-triangular preconditioners.
		\item We develop tailored two-grid solvers for the displacement and pressure blocks and use an inexpensive mass matrix to approximate the rotation block. Experiments in 2D and 3D assess parameter robustness and computational cost.
	\end{enumerate}

	The remainder of the paper is organized as follows. Section~\ref{sec:cv-biot} presents the Biot model, introduces the control-volume discretization
	and proves well-posedness of the fully discrete scheme. Section~\ref{sec:preconditioners} establishes the mechanics-block equivalence, analyzes the uncoupled block preconditioners, and develops the Schur-consistent block-diagonal and lower-triangular preconditioners together with their spectral and FOV estimates. Section~\ref{sec:numerics} reports the numerical experiments, and Section~\ref{sec:conclusion} gives the conclusions.

	\section{The control-volume discretization for the Biot model}
	\label{sec:cv-biot}

	\subsection{The mathematical model}

	Let \(d\in\{2,3\}\), and let \(\Omega\subset\mathbb R^d\) be a bounded polygonal domain for \(d=2\) or a bounded polyhedral domain for \(d=3\). We denote by \((\cdot,\cdot)_D\) the \(L^2(D)\) inner product and write \((\cdot,\cdot):=(\cdot,\cdot)_\Omega\). Boldface denotes vectors, and underlined symbols denote matrix-valued quantities. Let \(\underline\sigma_e\), \(\bs u\), \(\bs z\), and \(p\) denote the elastic stress, displacement, Darcy velocity, and pore pressure. The constitutive law is
	\begin{equation}
		\label{eq:elastic-constitutive}
		\Ac\underline\sigma_e=\varepsilon(\bs u),
		\qquad
		\varepsilon(\bs u):=\tfrac12(\nabla\bs u+\nabla\bs u^{\top}).
	\end{equation}
	For homogeneous isotropic elasticity, the compliance operator is
	\begin{equation}
		\label{eq:compliance}
		\Ac\underline\tau
		=\frac1{2\mu}
		\left(\underline\tau-
		\frac{\lambda}{2\mu+d\lambda}\tr(\underline\tau)\tI\right),
		\qquad
		\mu\ge\mu_{\min}>0,
		\quad \lambda\ge0.
	\end{equation}
	Here \(\lambda\) and \(\mu\) are the Lam\'e parameters describing volumetric and shear stiffness, respectively.
	The total poroelastic stress is
	\begin{equation}
		\label{eq:sigmap}
		\underline\sigma=\underline\sigma_e-\alpha_Bp\tI,
		\qquad 0<\alpha_B\le1.
	\end{equation}
	The Biot--Willis coefficient \(\alpha_B\) measures the coupling between pore pressure and stress.
	Given a body force \(\bs f\) and a fluid source \(q\), the quasi-static Biot system reads
	\begin{equation*}
		\begin{alignedat}{2}
			-\Div\underline\sigma&=\bs f
			&&\quad\text{in }\Omega\times(0,T],\\
			K^{-1}\bs z+\nabla p&=0
			&&\quad\text{in }\Omega\times(0,T],\\
			\partial_t\bigl(c_0p+\alpha_B\nabla\cdot\bs u\bigr)
			+\nabla\cdot\bs z&=q
			&&\quad\text{in }\Omega\times(0,T].
		\end{alignedat}
	\end{equation*}
	Here \(c_0\ge0\) is the storage coefficient, and the hydraulic mobility tensor \(K\) is symmetric and uniformly positive definite, with \(K_1|\bs\xi|^2\le K(x)\bs\xi\cdot\bs\xi\le K_2|\bs\xi|^2\) for a.e. \(x\in\Omega\) and all \(\bs\xi\in\mathbb R^d\), where \(0<K_1\le K_2\). We impose homogeneous displacement and pressure boundary conditions and prescribe \(p(0)=p^0\).

	To impose stress symmetry weakly, let \(m_d=d(d-1)/2\) and introduce the rotation variable \(\gamma\in\mathbb R^{m_d}\). For \(d=2\), \(\gamma=\tfrac12\trot(\bs u)\) and \(\asop(\underline\tau)=\tau_{21}-\tau_{12}\), whereas for \(d=3\), \(\gamma=\tfrac12\tcurl(\bs u)\) and \(\asop(\underline\tau)=(\tau_{32}-\tau_{23},\,\tau_{13}-\tau_{31},\,\tau_{21}-\tau_{12})^{\top}\). Using \eqref{eq:elastic-constitutive}--\eqref{eq:sigmap} and \(\nabla\cdot\bs u=\tr\bigl(\Ac(\underline\sigma+\alpha_Bp\tI)\bigr)\), we obtain the following weakly symmetric five-field formulation: Find \(\underline\sigma\in C^1([0,T];L^2(\Omega)^{d\times d})\), \(\bs u\in C([0,T];H_0^1(\Omega)^d)\), \(\gamma\in C([0,T];L^2(\Omega)^{m_d})\), \(\bs z\in C([0,T];L^2(\Omega)^d)\), and \(p\in C^1([0,T];H_0^1(\Omega))\) such that, for all \((\underline\psi,\bs v,\eta,\bs\xi,w)\in L^2(\Omega)^{d\times d}\times H_0^1(\Omega)^d\times L^2(\Omega)^{m_d}\times L^2(\Omega)^d\times H_0^1(\Omega)\),
	\begin{align*}
		(\Ac(\underline\sigma+\alpha_Bp\tI),\underline\psi)
		-(\nabla\bs u,\underline\psi)
		+(\gamma,\asop(\underline\psi))&=0,
		\\
		(\underline\sigma,\nabla\bs v)&=(\bs f,\bs v),
		\\
		(\asop(\underline\sigma),\eta)&=0,
		\\
		(K^{-1}\bs z,\bs\xi)+(\nabla p,\bs\xi)&=0,
		\\
		c_0(\partial_t p,w)
		+\alpha_B(\partial_t\Ac(\underline\sigma+\alpha_Bp\tI),w\tI)
		-(\bs z,\nabla w)&=(q,w).
	\end{align*}

	Unless stated otherwise, the hidden constants in \(\lesssim\), \(\gtrsim\), and \(\simeq\) are independent of the mesh size; dependence on physical parameters and on \(\Delta t\) is indicated explicitly.

	\subsection{Discrete spaces and fully discrete formulation}
	\label{subsec:five-field-cv}
Let \(\mathcal T_M\) be a shape-regular partition of \(\Omega\) into macro-cells. Each \(M\in\mathcal T_M\) is divided into \(2^d\) subcells by connecting an interior point of \(M\) to the midpoints of its edges in 2D or to the corresponding lower-dimensional face centers in 3D. The resulting subcell mesh is denoted by \(\mathcal T_h\). The split macro-facets form the set \(\mathcal F_{\rm pr}^{1/2}\), the facets internal to the macro-cells form \(\mathcal F_{\rm dl}\), and \(\mathcal F_h=\mathcal F_{\rm pr}^{1/2}\cup\mathcal F_{\rm dl}\). For each interior facet \(F=\partial E^+\cap\partial E^-\), let \(\normal\) point from \(E^+\) to \(E^-\) and write \(\jump{\phi}=\phi^+-\phi^-\); on boundary facets, \(\normal\) is the outward unit normal and the jump denotes the trace. 	{For \(M\in\mathcal T_M\), let \(\normal_M\) denote the outward unit normal on \(\partial M\).} The subcells sharing a primal vertex constitute an interaction region.
	For further details of the mesh construction, see \cite{FuZhao2025CiCP}.

	Define the normal-flux space
	\begin{equation*}
		\bs V_h
		:=\left\{
		\bs\xi_h\in\bs L^2(\Omega):
		\bs\xi_h|_E\in\bs P_0(E)\ \forall E\in\mathcal T_h,
		\quad
		\jump{\bs\xi_h\cdot\normal}_F=0\ \forall F\in
		\mathcal F_{\rm pr}^{1/2}\setminus\partial\Omega
		\right\}.
	\end{equation*}
	Its degrees of freedom are the oriented normal traces on the primal subfacets. Let \(P_0(\mathcal T_M):=\{w_h\in L^2(\Omega):w_h|_M\in P_0(M)\ \forall M\in\mathcal T_M\}\). The discrete spaces are \(\underline\Sigma_h=[\bs V_h]^d\), \(\bs Z_h=\bs V_h\), \(\bs U_h=[P_0(\mathcal T_M)]^d\), \(P_h=P_0(\mathcal T_M)\), and \(\Gamma_h=[P_0(\mathcal T_M)]^{m_d}\).

	For \(\underline\tau_h,\underline\psi_h\in\underline\Sigma_h\), \(\bs v_h\in \bs U_h\), \(\eta_h\in\Gamma_h\), \(\bs z_h,\bs\xi_h\in \bs Z_h\), and \(r_h,w_h\in P_h\), we define
	\begin{gather*}
		a_{\sigma\sigma}(\underline\tau_h,\underline\psi_h):=(\Ac\underline\tau_h,\underline\psi_h),\;
		a_{\sigma p}(r_h,\underline\psi_h):=(\alpha_B\Ac(r_h\tI),\underline\psi_h),\;
		a_{p\sigma}(\underline\tau_h,w_h):=a_{\sigma p}(w_h,\underline\tau_h),
		\\
		a_{pp}(r_h,w_h):=c_0(r_h,w_h)+\alpha_B^2(\Ac(r_h\tI),w_h\tI),\;
		a_{zz}(\bs z_h,\bs\xi_h):=(K^{-1}\bs z_h,\bs\xi_h).
	\end{gather*}
	The stress--displacement and Darcy--pressure couplings are
	\begin{align*}
		b_{\sigma u}(\bs v_h,\underline\psi_h)
		&:=\sum_{F\in\mathcal F_{\rm dl}}
		(\bs v_h,\jump{\underline\psi_h\normal})_F,
		&
		b_{u\sigma}(\underline\tau_h,\bs v_h)
		&:=\sum_{F\in\mathcal F_{\rm pr}^{1/2}}
		(\underline\tau_h\normal,\jump{\bs v_h})_F,
		\\
		b_{zp}(r_h,\bs\xi_h)
		&:=\sum_{F\in\mathcal F_{\rm dl}}
		(r_h,\jump{\bs\xi_h\cdot\normal})_F,
		&
		b_{pz}(\bs z_h,w_h)
		&:=\sum_{F\in\mathcal F_{\rm pr}^{1/2}}
		(\bs z_h\cdot\normal,\jump{w_h})_F.
	\end{align*}
	The weak-symmetry forms are \(c_{\sigma\gamma}(\eta_h,\underline\psi_h):=(\eta_h,\asop(\underline\psi_h))\) and \(c_{\gamma\sigma}(\underline\tau_h,\eta_h):=(\asop(\underline\tau_h),\eta_h)\).

    Let \(\mathcal F_M^{\rm int}\) denote the set of facets shared by two adjacent macrocells.
	For \(F=\partial M^+\cap\partial M^-\in\mathcal F_M^{\rm int}\),
	fix the unit normal \(\normal_F\) from \(M^+\) to \(M^-\) and define
	\(\jump{\phi}=\phi^+-\phi^-\) using the traces from \(M^\pm\).
    For any macro-facet or subfacet \(F\), set \(h_F:=\operatorname{diam}F\).
    In both 2D and 3D, we use full-jump rotation stabilization, penalizing all rotation components on every interior macro-facet. Define
	\begin{equation*}
	    j_\gamma(\eta_h,\xi_h):=
	    \sum_{F\in\mathcal F_M^{\rm int}}
	    h_F\bigl(\jump{\eta_h},\jump{\xi_h}\bigr)_F,
	    \qquad
	    |\eta_h|_{\gamma,h}^2:=j_\gamma(\eta_h,\eta_h),
	\end{equation*}
	and, for a fixed parameter \(\delta_\gamma>0\),
	set \(s_\gamma(\eta_h,\xi_h):=\delta_\gamma\mu\,j_\gamma(\eta_h,\xi_h)\). The stabilization is symmetric and positive semidefinite.
	Integration by parts gives
	\begin{equation}
		\label{eq:discrete-adjoint-identities}
		b_{\sigma u}(\bs v_h,\underline\tau_h)
		=-b_{u\sigma}(\underline\tau_h,\bs v_h),
		\;
        b_{zp}(w_h,\bs\xi_h)=-b_{pz}(\bs\xi_h,w_h),
		\;
		c_{\sigma\gamma}(\eta_h,\underline\tau_h)
		=c_{\gamma\sigma}(\underline\tau_h,\eta_h).
	\end{equation}

	Let \(t^n=n\Delta t\), and introduce the scaled Darcy variable \(\widehat{\bs z}_h^n:=\Delta t\,\bs z_h^n\). Given \((\underline\sigma_h^{n-1},p_h^{n-1})\), one backward Euler step seeks \( (\underline\sigma_h^n,\bs u_h^n,\gamma_h^n, \widehat{\bs z}_h^n,p_h^n) \in\underline\Sigma_h\times \bs U_h\times\Gamma_h\times \bs Z_h\times P_h \) such that, for all \((\underline\psi_h,\bs v_h,\eta_h,\bs\xi_h,w_h)\in\underline\Sigma_h\times \bs U_h\times\Gamma_h\times \bs Z_h\times P_h\),
	\begin{subequations}
		\begin{align}
			a_{\sigma\sigma}(\underline\sigma_h^n,\underline\psi_h)
			+a_{\sigma p}(p_h^n,\underline\psi_h)
			-b_{\sigma u}(\bs u_h^n,\underline\psi_h)
			+c_{\sigma\gamma}(\gamma_h^n,\underline\psi_h)&=0,
			\label{eq:be-discrete1}\\
			-b_{u\sigma}(\underline\sigma_h^n,\bs v_h)
			&=(\bs f^n,\bs v_h),
			\label{eq:be-discrete2}\\
			{c_{\gamma\sigma}(\underline\sigma_h^n,\eta_h)-s_\gamma(\gamma_h^n,\eta_h)}&{=0},
			\label{eq:be-discrete3}\\
			((\Delta t\,K)^{-1}\widehat{\bs z}_h^n,\bs\xi_h)
			+b_{zp}(p_h^n,\bs\xi_h)&=0,
			\label{eq:be-discrete4}\\
			a_{p\sigma}(\underline\sigma_h^n-\underline\sigma_h^{n-1},w_h)
			+a_{pp}(p_h^n-p_h^{n-1},w_h)
		+b_{pz}(\widehat{\bs z}_h^n,w_h)
			&=\Delta t\,(q^n,w_h).
			\label{eq:be-discrete5}
		\end{align}
	\end{subequations}
	Notably, the rotation-jump stabilization in the weak-symmetry equation \eqref{eq:be-discrete3} enables a mesh-uniform elasticity inf--sup estimate and yields a condensed rotation block uniformly spectrally equivalent to a \(\mu\)-scaled macrocell mass matrix.
	Taking \(\bs v_h=\bs e_i\chi_M\), \(i=1,\ldots,d\), and \(w_h=\chi_M\) in \eqref{eq:be-discrete2} and \eqref{eq:be-discrete5}, respectively, yields local momentum balance and fluid-mass conservation on each \(M\in\mathcal T_M\). Here \(\bs e_i\) is the \(i\)-th coordinate unit vector and \(\chi_M\) is the indicator function of \(M\).

	\subsection{Five-field and condensed three-field matrices}
	\label{subsec:matrix-forms}

	Let \(\{\underline\Psi_i\}\), \(\{\bs V_j\}\), \(\{\chi_k\}\), \(\{\bs\xi_\ell\}\), and \(\{\phi_m\}\) be bases of \(\underline\Sigma_h\), \(\bs U_h\), \(\Gamma_h\), \(\bs Z_h\), and \(P_h\), respectively. With the coefficient ordering \(\bm X=(\ssig,\uc,\gammac,\zc,\pp)^{\top}\), define
	\begin{equation}
		\label{eq:block-BC-0527}
		\begin{gathered}
			\begin{alignedat}{4}
				\A_{ij}&:=a_{\sigma\sigma}(\underline\Psi_j,\underline\Psi_i),&
				\quad\B_{ji}&:=b_{u\sigma}(\underline\Psi_i,\bs V_j),&
				\quad\C_{ki}&:=c_{\gamma\sigma}(\underline\Psi_i,\chi_k),&
				\quad\D_{m\ell}&:=b_{pz}(\bs\xi_\ell,\phi_m),
			\end{alignedat}
			\\
			\begin{alignedat}{3}
				\E_{im}&:=a_{\sigma p}(\phi_m,\underline\Psi_i),&
				\qquad\F_{\ell r}&:=((\Delta t\,K)^{-1}\bs\xi_r,\bs\xi_\ell),&
				\qquad(\App)_{mr}&:=\alpha_B^2(\Ac(\phi_r\tI),\phi_m\tI),
			\end{alignedat}
			\\
			\begin{alignedat}{3}
				(\Mzero)_{mr}&:=(\phi_r,\phi_m),&
				\qquad\Mp&:=c_0\Mzero,&
				\qquad\M&:=\Mp+\App.
			\end{alignedat}
		\end{gathered}
	\end{equation}
	\begingroup
	Let \(\Jg\) and \(\Mg\) denote the rotation stabilization and mass matrices, respectively:
	\begin{equation}
		\label{eq:Jgamma-Mgamma-matrices}
		(\Jg)_{k\ell}:=s_\gamma(\chi_\ell,\chi_k),
		\qquad
		(\Mg)_{k\ell}:=(\chi_\ell,\chi_k).
	\end{equation}
	\endgroup
	For isotropic elasticity, \(\App=\clambda\Mzero\), where \(\clambda:=\alpha_B^2/\Kdr\) and \(\Kdr:=\lambda+2\mu/d\). Here \(\Kdr\) is the drained bulk modulus. From \eqref{eq:discrete-adjoint-identities}, the five-field system takes the form
	\begin{equation}
		\label{eq:balanced-full-matrix-0527}
		\begin{bmatrix}
			\A&\B^{\top}&\C^{\top}&0&\E\\
			\B&0&0&0&0\\
			{\C}&0&{-\Jg}&0&0\\
			0&0&0&-\F&\D^{\top}\\
			\E^{\top}&0&0&\D&\M
		\end{bmatrix}
		\begin{bmatrix}
			\ssig\\ \uc\\ \gammac\\ \zc\\ \pp
		\end{bmatrix}
		=
		\begin{bmatrix}
			0\\-F_u\\0\\0\\r_p
		\end{bmatrix},
	\end{equation}
	where \((F_u)_j=(\bs f^n,\bs V_j)\) and \((r_p)_m=a_{p\sigma}(\underline\sigma_h^{n-1},\phi_m) +a_{pp}(p_h^{n-1},\phi_m)+\Delta t(q^n,\phi_m)\).  The diagonal block \(-\Jg\) is the matrix representation of \(-s_\gamma(\gamma_h,\eta_h)\).

Since \(\A\) and \(\F\) are symmetric positive definite and block diagonal over the interaction regions, their inverses are applied locally. Eliminating the stress and scaled Darcy flux gives
\begin{equation}
	\label{eq:stress-z-elimination}
	\ssig=-\A^{-1}(\B^{\top}\uc+\C^{\top}\gammac+\E\pp),
	\qquad
	\zc=\F^{-1}\D^{\top}\pp.
\end{equation}
Define the reduced blocks
\begin{equation}
	\label{eq:reduced-block-definitions}
	\begin{aligned}
		\Kuu&:=\B\A^{-1}\B^{\top},
		&
		\Kug&:=\B\A^{-1}\C^{\top},
		&
		\Kup&:=\B\A^{-1}\E,
		\\
		{\Kggzero}&{:=\C\A^{-1}\C^{\top}},
		&
		{\Kgg}&{:=\Kggzero+\Jg},
		&
		\Kgp&:=\C\A^{-1}\E,
		\\
		\Spsig&:=\E^{\top}\A^{-1}\E,
		&
		\Spz&:=\D\F^{-1}\D^{\top},
		&
		\Tp&:=\M-\Spsig+\Spz.
	\end{aligned}
\end{equation}
The pressure block \(\Tp\) is positive semidefinite. Indeed, \(\Mp=c_0\Mzero\succeq0\) and \(\Spz=\D\F^{-1}\D^{\top}\succeq0\). Moreover, \(\App-\Spsig\succeq0\), because it is the Schur complement with respect to \(\A\) of the block matrix \(\bigl[\begin{smallmatrix}\A&\E\\\E^{\top}&\App\end{smallmatrix}\bigr]\), whose quadratic form is \((\Ac(\underline\tau_h+\alpha_Bp_h\tI), \underline\tau_h+\alpha_Bp_h\tI)\). Consequently, \(\Tp=\Mp+(\App-\Spsig)+\Spz\succeq0\). Substituting the local elimination formulas \eqref{eq:stress-z-elimination} for \(\ssig\) and \(\zc\) into \eqref{eq:balanced-full-matrix-0527} gives the condensed displacement--rotation--pressure system
\begin{equation*}
	\begin{bmatrix}
		-\Kuu&-\Kug&-\Kup\\
		-\Kug^{\top}&-\Kgg&-\Kgp\\
		-\Kup^{\top}&-\Kgp^{\top}&\Tp
	\end{bmatrix}
	\begin{bmatrix}
		\uc\\ \gammac\\ \pp
	\end{bmatrix}
	=
	\begin{bmatrix}
		-F_u\\0\\r_p
	\end{bmatrix}.
\end{equation*}
Collecting the two mechanics variables in \(\mm=(\uc,\gammac)^{\top}\), set
\begin{equation}
	\label{eq:Km-Dm-L-def}
	\Kmzero:=
	\begin{bmatrix}
		\Kuu&\Kug\\
		\Kug^{\top}&\Kggzero
	\end{bmatrix},
	\qquad
	\Km:=\Kmzero+\diag(0,\Jg),
	\qquad
	\Dm:=\diag(\Kuu,\Kgg).
\end{equation}
Thus the condensed matrix is
\begin{equation}
	\label{eq:Ared-definition}
	\Ared=
	\begin{bmatrix}
		-\Km&-\Lop\\
		-\Lop^{\top}&\Tp
	\end{bmatrix}
	\qquad\text{with}\qquad
	\Lop=
	\begin{bmatrix}
		\Kup\\
		\Kgp
	\end{bmatrix}.
\end{equation}
\subsection{Well-posedness of the fully discrete scheme}
The following two mesh-uniform stability estimates are used to prove well-posedness of the fully discrete scheme. We use the following discrete norms:
\begin{equation}
	\label{eq:jump-norm}
	\|\varphi_h\|_{1,h}^2
	:=\sum_{F\in\mathcal F_{\rm pr}^{1/2}}
	h_F^{-1}\|\jump{\varphi_h}\|_{L^2(F)}^2,
	\qquad
	\|(\bs v_h,\eta_h)\|_h^2
	:=\|\bs v_h\|_{1,h}^2+\|\eta_h\|_{L^2(\Omega)}^2.
\end{equation}
Boundary jumps in \eqref{eq:jump-norm} incorporate the homogeneous Dirichlet conditions.

\begingroup
For the stability analysis below, consider a uniformly
shape-regular family of macrocell and subcell meshes.
In 2D, the subcells in each quadrilateral macrocell \(M\)
are \(\mathcal O(h_M^2)\) perturbations of
parallelograms~\cite{FuZhao2025CiCP}.
In 3D, we consider Cartesian macrocells, each divided into
eight corner cuboids by its coordinate midplanes.
Here \(h_M=\operatorname{diam}M\) and \(h=\max_M h_M\).

\begin{lemma}\relax
	\label{lem:method2-darcy-stability}
	\textup{(Darcy inf-sup condition).}
	There exists \(c_D>0\), depending only on the shape regularity of the mesh family, such that
	\begin{equation}
		\label{eq:method2-darcy-infsup}
		\sup_{\bm0\ne\bs\xi_h\in \bs Z_h}
		\frac{b_{zp}(w_h,\bs\xi_h)}{\|\bs\xi_h\|_{L^2(\Omega)}}
		\ge c_D\|w_h\|_{1,h}
		\qquad \forall w_h\in P_h.
	\end{equation}
\end{lemma}
The Darcy inf--sup estimate follows directly from \cite[eqs.~(4.14)--(4.15)]{FuZhao2025CiCP}, since
$\underline\Sigma_h=[\bs Z_h]^d$ and $\bs U_h=[P_h]^d$, we therefore omit the detailed proof.

\begin{lemma}[Elasticity inf-sup condition]
	\label{lem:method2-stabilized-elasticity-stability}
	There exist \(h_0>0\) and \(c_E>0\), independent of \(h\), \(\lambda\), and \(\mu\), such that, for all \(0<h\le h_0\),
	\begin{equation}
		\label{eq:method2-elasticity-stability}
		\left(
		\sup_{\underline0\ne\underline\tau_h\in\underline\Sigma_h}
		\frac{
			b_{u\sigma}(\underline\tau_h,\bs v_h)
			+c_{\gamma\sigma}(\underline\tau_h,\eta_h)
		}{\|\underline\tau_h\|_{L^2(\Omega)}}
		\right)^2
		+
		|\eta_h|_{\gamma,h}^2
		\ge c_E^2\|(\bs v_h,\eta_h)\|_h^2
	\end{equation}
	for all \((\bs v_h,\eta_h)\in \bs U_h\times\Gamma_h\).
\end{lemma}
\endgroup

\begingroup
\begin{proof}

	Set \(\mathcal K_h:=\{\underline\tau_h\in\underline\Sigma_h:
	b_{u\sigma}(\underline\tau_h,\bs w_h)=0, \forall\bs w_h\in \bs U_h\}\).
	In both dimensions we establish a stabilized divergence estimate and reconstruct a stress in \(\mathcal K_h\).  Invoking
	\cite[eqs.~(4.14)--(4.15)]{FuZhao2025CiCP}, we have \begin{equation}
		\label{eq:FuZhao-displacement-infsup}
		\|\bs v_h\|_{1,h}
		\le C
		\sup_{\underline0\ne\underline\tau_h\in\underline\Sigma_h}
		\frac{|b_{u\sigma}(\underline\tau_h,\bs v_h)|}
		{\|\underline\tau_h\|_{L^2(\Omega)}} \;\text{for any}\;\bs v_h\in \bs U_h.
	\end{equation}

	We next estimate the rotation component. We begin with 2D. For every $\eta_h\in\Gamma_h$, we prove
	\begin{equation}
		\label{eq:rotation-kernel-stabilized-estimate}
		\|\eta_h\|_{L^2(\Omega)}
		\le C\left(
		\sup_{\underline0\ne\underline\tau_h\in\mathcal K_h}
		\frac{|c_{\gamma\sigma}(\underline\tau_h,\eta_h)|}
		{\|\underline\tau_h\|_{L^2(\Omega)}}
		+|\eta_h|_{\gamma,h}
		\right).
	\end{equation}
	To prove this claim, let
	$\mathcal Q_h
	:=\{\bs\rho_h\in H_0^1(\Omega)^2:
	\bs\rho_h|_M\circ F_M\in[Q_1(\widehat M)]^2,
	\ \forall M\in\mathcal T_M\}$
	be the continuous bilinear quadrilateral finite element space on the macro-mesh, where \(F_M:\widehat M\to M\) is the
	bilinear element map from the reference square \(\widehat M\).
	We first show that, for \(q_h\in P_0(\mathcal T_M)\) satisfying
	\((q_h,1)=0\),
	\begin{equation}
		\label{eq:global-stabilized-Q1Q0}
		\|q_h\|_{L^2(\Omega)}
		\le C\left(
		\sup_{\bs0\ne\bs\rho_h\in\mathcal Q_h}
		\frac{|(\nabla\cdot\bs\rho_h,q_h)|}
		{\|\bs\rho_h\|_{H^1(\Omega)}}
		+|q_h|_{\gamma,h}
		\right).
	\end{equation}
	Choose \(\bs\rho\in H_0^1(\Omega)^2\) with
	\(\nabla\cdot\bs\rho=q_h\) and
	\(\|\bs\rho\|_{H^1(\Omega)}\le C\|q_h\|_{L^2(\Omega)}\), and let \(I_h\bs\rho\in\mathcal Q_h\)
	be a boundary-preserving Scott--Zhang-type quasi-interpolant. Such an
	operator is obtained from the local regularization construction for
	quadrilateral elements in \cite{BernardiGirault1998}, using the
	boundary-preserving treatment of \cite{ScottZhang1990}. The local
	\(H^1\)-stability and first-order \(L^2\)-approximation properties of
	this operator, together with the scaled trace inequality and the
	uniformly bounded overlap of element patches implies
	\(\sum_{F\in\mathcal F_M^{\rm int}}h_F^{-1}
	\|\bs\rho-I_h\bs\rho\|_{L^2(F)}^2\le C\|\bs\rho\|_{H^1(\Omega)}^2\).
	{Set \(\bs e=\bs\rho-I_h\bs\rho\).
	Since \(q_h\) is macrocellwise constant and \(\bs e|_{\partial\Omega}=0\),
	integration by parts and the Cauchy--Schwarz inequality give}
	\[
	\begin{aligned}
	\|q_h\|_{L^2(\Omega)}^2
	&\le |(q_h,\nabla\cdot I_h\bs\rho)|
	+\sum\nolimits_{F\in\mathcal F_M^{\rm int}}
	|(\jump{q_h},\bs e\cdot{\normal_F})_F|\\
	&\le |(q_h,\nabla\cdot I_h\bs\rho)|
	+|q_h|_{\gamma,h}
	\bigl(\sum\nolimits_{F\in\mathcal F_M^{\rm int}}
	h_F^{-1}\|\bs e\|_{L^2(F)}^2\bigr)^{1/2}\\
	&\le |(q_h,\nabla\cdot I_h\bs\rho)|
	+C|q_h|_{\gamma,h}\|q_h\|_{L^2(\Omega)},
	\end{aligned}
	\]
    which proves \eqref{eq:global-stabilized-Q1Q0}.

	We next prove \eqref{eq:rotation-kernel-stabilized-estimate} using \eqref{eq:global-stabilized-Q1Q0}. Write $\nabla^\perp r:=(\partial_y r,-\partial_x r)$ for a scalar $r$.
	For $\bs\rho_h\in\mathcal Q_h$, define
	$\mathcal R_h\bs\rho_h\in\underline\Sigma_h$
	through its row-wise normal traces on the primal subfacets:
	\begin{equation}
		\label{eq:FuZhao-rotated-gradient-transfer}
		\int_F (\mathcal R_h\bs\rho_h)_i\cdot\normal\,ds
		=
		\int_F \nabla^\perp(\rho_h)_i\cdot\normal\,ds,
		\qquad i=1,2,
		\quad F\in\mathcal F_{\rm pr}^{1/2}.
	\end{equation}
	We use one common value of the flux degree of freedom on each interior primal half-edge, with opposite orientations on the two adjacent subcells. Thus \eqref{eq:FuZhao-rotated-gradient-transfer} defines
	a member of \(\underline\Sigma_h\) directly in physical coordinates.
	On the other hand, standard scaling gives $\|\mathcal R_h\bs\rho_h\|_{L^2(\Omega)}\le C\|\bs\rho_h\|_{H^1(\Omega)}$.
	Moreover, for any $M\in\mathcal T_M$ and $i=1,2$, we have $\int_{\partial M}(\mathcal R_h\bs\rho_h)_i\cdot{\normal_M}\,ds=0$.
	Since every $\bs w_h\in \bs U_h$ is constant on $M$, summation over the macrocells yields
	$b_{u\sigma}(\mathcal R_h\bs\rho_h,\bs w_h)=0$, so $\mathcal R_h\bs\rho_h\in\mathcal K_h$.

	On a parallelogram \(M\) with the center-to-midpoint subdivision,
	let \(K_v\) be the corner subcell at vertex \(v\).
	Since the edge traces of \((\rho_h)_i\) are affine,
	\eqref{eq:FuZhao-rotated-gradient-transfer} matches the normal
	components of \((\mathcal R_h\bs\rho_h)_i|_{K_v}\) and
	\(\nabla^\perp(\rho_h)_i(v)\) on the two adjacent primal half-edges.
	Their normals are linearly independent, so
	\((\mathcal R_h\bs\rho_h)_i|_{K_v}
	=\nabla^\perp(\rho_h)_i(v)\), \(i=1,2\).
	Thus
	\(\asop((\mathcal R_h\bs\rho_h)|_{K_v})
	=\nabla\cdot\bs\rho_h(v)\).
	The divergence of a bilinear \(Q_1\) function, pulled back by the
	affine map from the reference square, belongs to \(Q_1\). Since the four
	subcells have equal area, the tensor-product trapezoidal rule is exact
	for this function, and hence
	\begin{equation}
		\label{eq:FuZhao-transfer-exact-parallelogram1}
		\int_M\asop(\mathcal R_h\bs\rho_h)\,dx
		=
		\int_M\nabla\cdot\bs\rho_h\,dx.
	\end{equation}
	For a perturbed quadrilateral \(M\), the preceding quadrature argument
	is not exact. Let \(M^0\) be its parallelogram counterpart, with
	\(|v-v^0|\le Ch_M^2\) at corresponding vertices, and let
	\(\bs\rho_h^0\) have the same nodal values as \(\bs\rho_h\).
	Subtracting the common value at one vertex \(v_*\) and applying
	finite-dimensional scaling gives
	\begin{equation}
		\label{eq:affine-counterpart-scaling}
		|\bs\rho_h^0|_{H^1(M^0)}
		+\max_v|\bs\rho_h(v)-\bs\rho_h(v_*)|
		\le C|\bs\rho_h|_{H^1(M)}.
	\end{equation}
	On a corner subcell \(K_v\), let \(\underline{N}_v\) contain the integrated normals
	of its two primal half-edges and \(\bs d_{v,i}\) their rotated-gradient
	fluxes. Use superscript \(0\) for \(M^0\), and write \(\bs r_{v,i}\) and
	\(\bs r_{v,i}^0\) for the constant reconstructed rows. The affine edge
	traces have matching endpoint values, so \(\bs d_{v,i}=\bs d_{v,i}^0\) and
	\(\underline{N}_v^{\top}\bs r_{v,i}=\bs d_{v,i}=(\underline{N}_v^0)^{\top}\bs r_{v,i}^0\). Uniform shape regularity,
	\eqref{eq:affine-counterpart-scaling}, and the vertex perturbation give
	\[
	|\bs d_{v,i}|\le C|\bs\rho_h|_{H^1(M)},\quad
	\|\underline{N}_v-\underline{N}_v^0\|_2\le Ch_M^2,\quad
	\|(\underline{N}_v^{\top})^{-1}\|_2+\|((\underline{N}_v^0)^{\top})^{-1}\|_2\le Ch_M^{-1},
	\]
where \(\|\cdot\|_2\) is the induced Euclidean matrix norm. Thus \(|\bs r_{v,i}-\bs r_{v,i}^0|\le C|\bs\rho_h|_{H^1(M)}\) and
	\(\bigl||K_v|-|K_v^0|\bigr|\le Ch_M^3\).
	Because the reconstructed stresses are constant on the corner subcells,
	their integrals are sums of subcell areas times the corresponding
	asymmetries. Using \(|\bs r_{v,i}^0|\le Ch_M^{-1}|\bs\rho_h|_{H^1(M)}\),
	the row difference and area bounds above give
	$\left|\int_M\asop(\mathcal R_h\bs\rho_h)\,dx
	-\int_{M^0}\asop(\mathcal R_h^0\bs\rho_h^0)\,dx\right|
	\le C h_M^2|\bs\rho_h|_{H^1(M)}$.
	The divergence theorem gives the same bound for the difference of the
	divergence integrals: subtract the common nodal value, use
	\eqref{eq:affine-counterpart-scaling} for the remaining affine edge traces,
	and note that corresponding integrated edge normals differ by \({\mathcal O}(h_M^2)\).
	Combining the two bounds with
	\eqref{eq:FuZhao-transfer-exact-parallelogram1} on \(M^0\) yields
	\begin{equation}
		\label{eq:FuZhao-transfer-perturbation-local}
		\left|\int_M\bigl(\asop(\mathcal R_h\bs\rho_h)
		-\nabla\cdot\bs\rho_h\bigr)\,dx\right|
		\le C h_M^2|\bs\rho_h|_{H^1(M)}.
	\end{equation}

    	 Since \(q_h\) is constant on each macrocell, summation of
	\eqref{eq:FuZhao-transfer-perturbation-local} and the
	Cauchy--Schwarz inequality yield the global defect explicitly as follows.
	Writing \(q_M=q_h|_M\), using
	\(\|q_h\|_{L^2(M)}\simeq h_M|q_M|\), we obtain
    \begin{equation}
        \label{eq:FuZhao-transfer-perturbation-global}
        \begin{aligned}
            &\left|c_{\gamma\sigma}(\mathcal R_h\bs\rho_h,q_h)
            -(\nabla\cdot\bs\rho_h,q_h)\right|
            \le C\sum_{M\in\mathcal T_M}|q_M|h_M^2|\bs\rho_h|_{H^1(M)}\\
            &\qquad\le Ch\sum_{M\in\mathcal T_M}
            \|q_h\|_{L^2(M)}|\bs\rho_h|_{H^1(M)}
            \le Ch\|q_h\|_{L^2(\Omega)}\|\bs\rho_h\|_{H^1(\Omega)}.
        \end{aligned}
    \end{equation}

	Write \(\eta_h=\eta_h^0+\bar\eta_h\), where
	\(\bar\eta_h=|\Omega|^{-1}\int_\Omega\eta_h\).
	Set \(\underline{S}(a)=\frac a2\bigl(\begin{smallmatrix}0&-1\\1&0\end{smallmatrix}\bigr)\)
	and \(m_h=|\Omega|^{-1}\int_\Omega\asop(\mathcal R_h\bs\rho_h)\).
	The constant stress \(\underline{S}(m_h)\in\mathcal K_h\), since its
	normal flux sums to zero on every macrocell.
	Thus \(\mathcal Z_h\bs\rho_h:=
	\mathcal R_h\bs\rho_h-\underline{S}(m_h)\in\mathcal K_h\), with
	\(\int_\Omega\asop(\mathcal Z_h\bs\rho_h)\,dx=0\) and
	\mbox{\(\|\mathcal Z_h\bs\rho_h\|_{L^2(\Omega)}
	\le C\|\bs\rho_h\|_{H^1(\Omega)}\)}.
	Since \(\asop \underline{S}(m_h)=m_h\) and
	\(\int_\Omega\eta_h^0=0\),
	$
	c_{\gamma\sigma}(\mathcal Z_h\bs\rho_h,\eta_h)
	=c_{\gamma\sigma}(\mathcal R_h\bs\rho_h,\eta_h^0).
	$
	Define
	$
	\mathfrak S_h(\eta_h):=
	\sup_{\underline0\ne\underline\tau_h\in\mathcal K_h}
	\frac{|c_{\gamma\sigma}(\underline\tau_h,\eta_h)|}
	{\|\underline\tau_h\|_{L^2(\Omega)}}.
	$
	Taking \(q_h=\eta_h^0\) in \eqref{eq:FuZhao-transfer-perturbation-global} and using
	\(\|\mathcal Z_h\bs\rho_h\|_{L^2{(\Omega)}}\le C\|\bs\rho_h\|_{H^1{(\Omega)}}\), we obtain, for every
	\(\bs\rho_h\in\mathcal Q_h\),
\[
	\begin{aligned}
		|(\nabla\cdot\bs\rho_h,\eta_h^0)|
		&\le |c_{\gamma\sigma}(\mathcal Z_h\bs\rho_h,\eta_h)|
		+Ch\|\bs\rho_h\|_{H^1{(\Omega)}}\|\eta_h^0\|_{L^2{(\Omega)}}\\
		&\le C\bigl(\mathfrak S_h(\eta_h)
		+h\|\eta_h^0\|_{L^2{(\Omega)}}\bigr)
		\|\bs\rho_h\|_{H^1{(\Omega)}}.
	\end{aligned}
	\]
	Moreover,
	\(|\eta_h^0|_{\gamma,h}=|\eta_h|_{\gamma,h}\). Hence
	\eqref{eq:global-stabilized-Q1Q0} implies
    $\|\eta_h^0\|_{L^2{(\Omega)}}
	\le C\bigl(\mathfrak S_h(\eta_h)+|\eta_h|_{\gamma,h}\bigr)
	+C_{\rm def}h\|\eta_h^0\|_{L^2{(\Omega)}}$,
	where \(C_{\rm def}\) is independent of \(h\). Choose \(h_0>0\) so that
	\(C_{\rm def}h_0\le\tfrac12\). For \(h\le h_0\), moving the last term
	to the left gives
	$\|\eta_h^0\|_{L^2(\Omega)}
	\le C\bigl(\mathfrak S_h(\eta_h)+|\eta_h|_{\gamma,h}\bigr)$.

	{We next bound \(\bar\eta_h\).} Since
	\(\asop \underline{S}(1)=1\), \(\|\underline{S}(1)\|_{L^2(\Omega)}^2=|\Omega|/2\), and
	$
	c_{\gamma\sigma}(\underline{S}(1),\eta_h)
	=\int_\Omega\eta_h\,dx
	=|\Omega|\bar\eta_h,
	$
	we have
	$
	\|\bar\eta_h\|_{L^2(\Omega)}
	=|\Omega|^{1/2}|\bar\eta_h|
	=\frac1{\sqrt2}
	\frac{|c_{\gamma\sigma}(\underline{S}(1),\eta_h)|}
	{\|\underline{S}(1)\|_{L^2(\Omega)}}
	\le \frac1{\sqrt2}\mathfrak S_h(\eta_h).
	$
	Combining the bounds for \(\eta_h^0\) and \(\bar\eta_h\)
proves \eqref{eq:rotation-kernel-stabilized-estimate}.

	{By the definition of \(\mathcal K_h\), \eqref{eq:rotation-kernel-stabilized-estimate} gives}
	\begin{equation}
		\label{eq:eta-by-coupling-and-jump}
		\|\eta_h\|_{L^2(\Omega)}
		\le C\left(
		\sup_{\underline0\ne\underline\tau_h\in\underline\Sigma_h}
		\frac{|b_{u\sigma}(\underline\tau_h,\bs v_h)+c_{\gamma\sigma}(\underline\tau_h,\eta_h)|}
		{\|\underline\tau_h\|_{L^2(\Omega)}}
		+|\eta_h|_{\gamma,h}
		\right).
	\end{equation}
	By \eqref{eq:FuZhao-displacement-infsup} and the Cauchy--Schwarz inequality,
	\(\|\bs v_h\|_{1,h}\) is controlled by the coupling supremum in
	\eqref{eq:eta-by-coupling-and-jump} and
	\(\|\eta_h\|_{L^2(\Omega)}\). Combining this with
	\eqref{eq:eta-by-coupling-and-jump} proves
	\eqref{eq:method2-elasticity-stability} in 2D.

	For \(d=3\), let \(\mathcal Q_h^{(3)}\subset H_0^1(\Omega)\) be the continuous trilinear \(Q_1\) space on the Cartesian macro-mesh, and let \(\mathcal Q_h^{(3),3\times3}\) be its matrix-valued counterpart. {For matrix-valued fields, both \(\Div\) and \(\tcurl\) act row by row.} We first establish the vector-valued stabilized divergence estimate, for \(\bs q_h\in\Gamma_h\) with zero mean componentwise,
	\begin{equation}
		\label{eq:3d-stabilized-Q1Q0}
		\|\bs q_h\|_{L^2(\Omega)}
		\le C\left(
		\sup_{{\underline0}\ne \underline{W}_h\in\mathcal Q_h^{(3),3\times3}}
		\frac{|(\bs q_h,\operatorname{div}\underline{W}_h)|}{\|\underline{W}_h\|_{H^1(\Omega)}}
		+|\bs q_h|_{\gamma,h}\right).
	\end{equation}
	Choose \(\underline{W}\in H_0^1(\Omega)^{3\times3}\) row by row so that
	\(\operatorname{div}\underline{W}=\bs q_h\) and \(\|\underline{W}\|_{H^1(\Omega)}\le C\|\bs q_h\|_{L^2(\Omega)}\).
	The boundary-preserving quasi-interpolation construction of
	\cite[Section~6]{ErnGuermond2017}, applied componentwise and specialized
	to the conforming tensor-product \(Q_1\) element on affine Cartesian
	cuboids, provides \(I_h\underline{W}\in\mathcal Q_h^{(3),3\times3}\) satisfying
	$\|I_h\underline{W}\|_{H^1(\Omega)}^2+
	\sum_{F\in\mathcal F_M^{\rm int}}h_F^{-1}
	\|\underline{W}-I_h\underline{W}\|_{L^2(F)}^2\le C\|\underline{W}\|_{H^1(\Omega)}^2$.
	Indeed, with \(\underline{e}=\underline{W}-I_h\underline{W}\) and a fixed-size patch \(\omega_M\), the local estimates
	\(\|\underline{e}\|_{L^2(M)}\le Ch_M\|\underline{W}\|_{H^1(\omega_M)}\) and
	\(\|\nabla \underline{e}\|_{L^2(M)}\le C\|\underline{W}\|_{H^1(\omega_M)}\) and the scaled trace inequality give, on every face \(F\subset\partial M\),
	$h_F^{-1}\|\underline{e}\|_{L^2(F)}^2
	\le C\bigl(h_M^{-2}\|\underline{e}\|_{L^2(M)}^2
	+\|\nabla \underline{e}\|_{L^2(M)}^2\bigr)
	\le C\|\underline{W}\|_{H^1(\omega_M)}^2$.
	Summation proves the stated face bound because the patches have uniformly bounded overlap. Since \(\bs q_h\) is macrocellwise constant and \(\underline{W}-I_h\underline{W}\) has zero boundary trace,
	{macrocell integration by parts gives
	\((\bs q_h,\Div \underline{e})=\sum_{F\in\mathcal F_M^{\rm int}}
	(\jump{\bs q_h},\underline{e}\normal_F)_F\).
	The Cauchy--Schwarz inequality then yields}
	\[
	\|\bs q_h\|_{L^2(\Omega)}^2
	\le |(\bs q_h,\operatorname{div}I_h\underline{W})|
	+|\bs q_h|_{\gamma,h}
	\left(\sum_{F\in\mathcal F_M^{\rm int}}h_F^{-1}
	\|\underline{W}-I_h\underline{W}\|_{L^2(F)}^2\right)^{1/2}.
	\]
	This proves \eqref{eq:3d-stabilized-Q1Q0}.

	Define the linear map
	\(\Xi\underline{\Phi}:=\underline{\Phi}^{\top}
	-(\operatorname{tr}\underline{\Phi}){\tI}\).
	For \(\underline{W}_h\in\mathcal Q_h^{(3),3\times3}\), set
	$\underline{\Phi}_h:=\Xi^{-1}\underline{W}_h
	=\underline{W}_h^{\top}
	-\tfrac12(\operatorname{tr}\underline{W}_h){\tI}$.
Using the definition of \(\asop\) and the row-wise curl, we obtain
	\[
	\bigl(\asop(\operatorname{curl}\underline{\Phi}_h)\bigr)_i
	=\sum_{j=1}^3\partial_j(\underline{\Phi}_h)_{ji}
	-\partial_i\operatorname{tr}\underline{\Phi}_h
	=\bigl(\operatorname{div}(\Xi\underline{\Phi}_h)\bigr)_i,
	\qquad i=1,2,3.
	\]
	Hence \(\asop(\operatorname{curl}\underline{\Phi}_h)=\operatorname{div}\underline{W}_h\).

	{For each corner cuboid \(E\), let \(F_j(E)\) be its primal quarterface normal to direction \(j\).
	We then define \(\mathcal R_h\underline{\Phi}_h\in\underline\Sigma_h\) to be constant on each \(E\), with
	\((\mathcal R_h\underline{\Phi}_h)_{ij}|_E:=|F_j(E)|^{-1}\int_{F_j(E)}(\operatorname{curl}\underline{\Phi}_h)_{ij}\,ds\), \(i,j=1,2,3\).
	This uniquely determines the reconstructed stress and preserves the row-wise normal flux integrals on every primal quarterface.}
	Across an interior primal quarterface, the normal curl flux is single-valued: it depends only on tangential derivatives of the common trace of \(\underline{\Phi}_h\). The reconstructed rows therefore have matching normal components. {Let \(h_E=\operatorname{diam}E\).} Then
	\[
	\|{\mathcal R_h}\underline{\Phi}_h\|_{L^2(E)}^2
	\le C h_E\sum_{i,j=1}^3
	\|(\operatorname{curl}\underline{\Phi}_h)_{ij}\|_{L^2(F_j(E))}^2
	\le C\|\nabla\underline{\Phi}_h\|_{L^2(E)}^2.
	\]
	Summing over \(E\), and using that \(\Xi^{-1}\) is fixed and bounded, yields
	\begin{equation}
		\label{eq:3d-curl-transfer-bound}
		\|{\mathcal R_h}\underline{\Phi}_h\|_{L^2(\Omega)}\le C\|\underline{W}_h\|_{H^1(\Omega)}.
	\end{equation}
	{By the definition of \(\mathcal R_h\), we have
	\(\int_{\partial M}(\mathcal R_h\underline{\Phi}_h)_i\cdot\normal_M\,ds
	=\int_{\partial M}(\operatorname{curl}\underline{\Phi}_h)_i\cdot\normal_M\,ds=0\)
	for every \(M\in\mathcal T_M\) and \(i=1,2,3\). Therefore \(\mathcal R_h\underline{\Phi}_h\in\mathcal K_h\).}

	We next show that the reconstruction preserves the integral of the asymmetry on each Cartesian macrocell.
	Write \(g_{ij}=(\operatorname{curl}\underline{\Phi}_h)_{ij}\). Let \(F_j^-\) and \(F_j^+\) be the opposite macrofaces normal to direction \(j\), separated by the side length \(\ell_j\). Each is partitioned into four primal quarterfaces. On a corner cuboid \(E\) adjoining a quarterface \(F\) normal to direction \(j\), flux matching and \(|E|=(\ell_j/2)|F|\) give
	$
	\int_E({\mathcal R_h}\underline{\Phi}_h)_{ij}\,dx
	=\frac{\ell_j}{2}\int_F g_{ij}\,ds.
	$
	Summing over the eight corner cuboids gives
	\begin{equation}
		\label{eq:3d-exact-curl-transfer}
		\int_M({\mathcal R_h}\underline{\Phi}_h)_{ij}\,dx
		=\frac{\ell_j}{2}\left(
		\int_{F_j^-}g_{ij}\,ds+\int_{F_j^+}g_{ij}\,ds\right)
		=\int_Mg_{ij}\,dx.
	\end{equation}
	The final equality is the 1D trapezoidal identity in \(x_j\): trilinearity of \(\underline{\Phi}_h\) makes \(g_{ij}\) affine in that coordinate, although it may vary tangentially.
	Because rotations are constant on each macrocell,
	\eqref{eq:3d-exact-curl-transfer} and
	\(\asop(\operatorname{curl}\underline{\Phi}_h)
	=\operatorname{div}\underline{W}_h\) imply
	\begin{equation}
		\label{eq:3d-curl-coupling-identity}
		{c_{\gamma\sigma}(\mathcal R_h\underline{\Phi}_h,\eta_h)
		=(\operatorname{div}\underline{W}_h,\eta_h),
		\qquad\forall\eta_h\in\Gamma_h.}
	\end{equation}
	Because \(\underline{W}_h\) has zero boundary trace, \((\operatorname{div}\underline{W}_h,\bar\eta_h)=0\), where \(\bar\eta_h=|\Omega|^{-1}\int_\Omega\eta_h\). Applying \eqref{eq:3d-stabilized-Q1Q0} to \(\eta_h^0=\eta_h-\bar\eta_h\), and then using {\eqref{eq:3d-curl-coupling-identity}} and \eqref{eq:3d-curl-transfer-bound}, gives
	$
	\|\eta_h^0\|_{L^2(\Omega)}
	\le C\left(
	{\mathfrak S_h(\eta_h)}
	+|\eta_h|_{\gamma,h}\right).
	$
	{We next bound \(\bar\eta_h\).} For \(\bs a\in\mathbb R^3\), the globally constant skew stress
	${\underline{S}(\bs a)=\tfrac12\bigl(\begin{smallmatrix}0&-a_3&a_2\\a_3&0&-a_1\\-a_2&a_1&0\end{smallmatrix}\bigr)}$
	belongs to \(\mathcal K_h\) and satisfies \(\asop \underline{S}(\bs a)=\bs a\). {Choosing \(\underline\tau_h=\underline{S}(\bar\eta_h)\) in the definition of \(\mathfrak S_h(\eta_h)\) gives \(\|\bar\eta_h\|_{L^2(\Omega)}\le C\mathfrak S_h(\eta_h)\).} Consequently \eqref{eq:rotation-kernel-stabilized-estimate} holds also in 3D.

	{Finally, for \(\bs v_h\in \bs U_h\), prescribe \(\underline\tau_{v,h}\normal=h_F^{-1}\jump{\bs v_h}\) on each primal quarterface \(F\), with boundary jumps interpreted as traces.
	The three linearly independent primal-quarterface normals uniquely determine a constant stress on each corner cuboid.
	The prescribed normal traces agree across interior primal quarterfaces, hence \(\underline\tau_{v,h}\in\underline\Sigma_h\).} Scaling gives
 	$b_{u\sigma}(\underline\tau_{v,h},\bs v_h)=\|\bs v_h\|_{1,h}^2$ and
	$\|\underline\tau_{v,h}\|_{L^2(\Omega)}\le C\|\bs v_h\|_{1,h}$.
	This is the 3D displacement lift.
	Combining this displacement estimate with the preceding bounds for \(\eta_h^0\) and \(\bar\eta_h\) proves \eqref{eq:method2-elasticity-stability}.
	\unskip\nobreak\hfill\proofSymbol\NoEndMark
\end{proof}

\endgroup

\begin{theorem}\relax
	\label{thm:fully-discrete-wellposedness}
	\textup{(Well-posedness).}
	The system \eqref{eq:be-discrete1}--\eqref{eq:be-discrete5} has a unique solution.
\end{theorem}

\begin{proof}
	The system is square and finite dimensional, so it suffices to prove that its homogeneous problem has only the zero solution. Let \((\underline\sigma_h^*,\bs u_h^*,\gamma_h^*,\widehat{\bs z}_h^*,p_h^*)\in\underline\Sigma_h\times \bs U_h\times\Gamma_h\times \bs Z_h\times P_h\) satisfy the homogeneous system. \begingroup Taking \(\underline\psi_h=\underline\sigma_h^*\), \(\bs v_h=\bs u_h^*\), and \(\eta_h=\gamma_h^*\) in \eqref{eq:be-discrete1}--\eqref{eq:be-discrete3}, respectively, and using \eqref{eq:discrete-adjoint-identities}, we obtain
	\begin{equation}
		\label{eq:wellposed-mechanical-kernel}
		a_{\sigma\sigma}(\underline\sigma_h^*,\underline\sigma_h^*)+a_{\sigma p}(p_h^*,\underline\sigma_h^*)+s_\gamma(\gamma_h^*,\gamma_h^*)=0.
	\end{equation}
	\endgroup
	Likewise, taking \(\bs\xi_h=\widehat{\bs z}_h^*\) and \(w_h=p_h^*\) in \eqref{eq:be-discrete4} and \eqref{eq:be-discrete5}, and adding the two equalities, yields
	\begin{equation}
		\label{eq:wellposed-flow-kernel}
		((\Delta t\,K)^{-1}\widehat{\bs z}_h^*,\widehat{\bs z}_h^*)+a_{p\sigma}(\underline\sigma_h^*,p_h^*)+a_{pp}(p_h^*,p_h^*)=0.
	\end{equation}
	Adding \eqref{eq:wellposed-mechanical-kernel} and \eqref{eq:wellposed-flow-kernel} yields
	\begin{equation*}
		0=\|\Ac^{1/2}(\underline\sigma_h^*+\alpha_Bp_h^*\tI)\|_{L^2(\Omega)}^2+c_0\|p_h^*\|_{L^2(\Omega)}^2+((\Delta t\,K)^{-1}\widehat{\bs z}_h^*,\widehat{\bs z}_h^*)+s_\gamma(\gamma_h^*,\gamma_h^*).
	\end{equation*}
	By the positive definiteness of \(\Ac\) and the nonnegativity of the stabilization, we obtain
	\begin{equation}
		\label{eq:wellposed-first-zero}
		\underline\sigma_h^*+\alpha_Bp_h^*\tI=\underline0,\qquad \widehat{\bs z}_h^*=\bs0,\qquad |\gamma_h^*|_{\gamma,h}=0.
	\end{equation}
	Together, \eqref{eq:be-discrete4} and \eqref{eq:wellposed-first-zero} imply \(b_{zp}(p_h^*,\bs\xi_h)=0\) for all \(\bs\xi_h\in \bs Z_h\). Hence, by \eqref{eq:method2-darcy-infsup},
	\[
	c_D\|p_h^*\|_{1,h}\le\sup_{\bm0\ne\bs\xi_h\in \bs Z_h}\frac{b_{zp}(p_h^*,\bs\xi_h)}{\|\bs\xi_h\|_{L^2(\Omega)}}=0.
	\]
	The discrete Poincar\'e inequality for \(P_0(\mathcal T_M)\) gives \(\|p_h^*\|_{L^2(\Omega)}\lesssim \|p_h^*\|_{1,h}=0\). Thus \(p_h^*=0\), and \eqref{eq:wellposed-first-zero} yields \(\underline\sigma_h^*=\underline0\).

	Finally, \eqref{eq:be-discrete1} reduces to \(b_{u\sigma}(\underline\psi_h,\bs u_h^*)+c_{\gamma\sigma}(\underline\psi_h,\gamma_h^*)=0\) for all \(\underline\psi_h\in\underline\Sigma_h\). Together with \(|\gamma_h^*|_{\gamma,h}=0\), the stabilized condition \eqref{eq:method2-elasticity-stability} yields \(\|(\bs u_h^*,\gamma_h^*)\|_h=0\). Hence \((\underline\sigma_h^*,\bs u_h^*,\gamma_h^*,\widehat{\bs z}_h^*,p_h^*)=\bm0\).
\end{proof}

\section{Block preconditioners and spectral analysis}
\label{sec:preconditioners}
This section develops block preconditioners for the reduced displacement--rotation--pressure matrix \eqref{eq:Ared-definition}. Starting from the parameter-robust equivalence \(\Km\simeq\Dm\), {which} justifies replacing the coupled mechanics block by separate displacement and rotation blocks, we study the inexpensive uncoupled block-diagonal preconditioners
\(\mathcal P_{\rm bd}^{\Tp}\) and \(\mathcal P_{\rm bd}^{\Pmz}\). They retain the natural pressure blocks, and their stability is characterized by the exact discrete ratio \(\Thetah\), with the explicit sufficient estimate \(\Thetah\le \alpha_B^2/(c_0\Kdr)\) when \(c_0>0\).
The third subsection introduces the exact and Schur-consistent pressure blocks. We prove uniform spectral equivalence for the block-diagonal preconditioners \(\mathcal P_{\rm bd}^{\rm ex}\) and \(\mathcal P_{\rm bd}^{\rm sc}\). Moreover, we derive weighted FOV estimates for the lower-triangular factors \(\mathcal P_{\rm lt}^{\rm ex}\) and \(\mathcal P_{\rm lt}^{\rm sc}\).
The distinction is cost versus robustness: the uncoupled preconditioners are cheaper but require pressure control through \(\Tp\), whereas the Schur-consistent preconditioners use a more expensive pressure block.
Finally, we describe practical implementations of the displacement, rotation, and pressure inverses.

\subsection{Mechanical block equivalence}
The purpose of this subsection is to prove
\[
\Km\simeq\Dm,
\]
with constants independent of \(h\), \(\lambda\), and \(\mu\). The upper bound is immediate once the condensed energies are written in terms of compliance-Riesz representatives. The essential point is the lower bound: the displacement and rotation representatives may in principle cancel in the coupled energy, while a direct uniform comparison of the compliance and \(L^2\) stress norms is unavailable as \(\lambda\to\infty\). The inf--sup estimate controls the rotation through the total stress representative and the rotation-jump stabilization term.

We use the compliance norm \(\norm{\underline\tau_h}_{\Ac}^2:=a_{\sigma\sigma}(\underline\tau_h,\underline\tau_h)\). For \(\underline\tau_h\in\underline\Sigma_h\), let \(\underline\tau_h^D:=\underline\tau_h-\frac1d\tr(\underline\tau_h)\tI\) denote its deviatoric part. \begingroup The \(L^2\)-orthogonal decomposition into deviatoric and spherical parts gives
$\norm{\underline\tau_h}_{L^2(\Omega)}^2=\norm{\underline\tau_h^D}_{L^2(\Omega)}^2+\frac1d\norm{\tr(\underline\tau_h)}_{L^2(\Omega)}^2$.
Combining this identity with the isotropic compliance tensor \eqref{eq:compliance}, we obtain
\[
\norm{\underline\tau_h}_{\Ac}^2=\frac1{2\mu}\norm{\underline\tau_h^D}_{L^2(\Omega)}^2+\frac1{d(2\mu+d\lambda)}\norm{\tr(\underline\tau_h)}_{L^2(\Omega)}^2.
\]
Consequently, 	$\norm{\underline\tau_h^D}_{L^2(\Omega)}\le\sqrt{2\mu}\,\norm{\underline\tau_h}_{\Ac}$, and	$\norm{\underline\tau_h}_{\Ac}\le\frac1{\sqrt{2\mu}}\norm{\underline\tau_h}_{L^2(\Omega)}$.
\endgroup

{For any coefficient vector \(\mm=(\uc,\gammac)\), set
\(\uu:=\sum_j(\uc)_j\,\bs V_j\) and
\(\gh:=\sum_k(\gammac)_k\,\chi_k\)
using the bases \(\{\bs V_j\}\) of \(\bs U_h\) and \(\{\chi_k\}\) of \(\Gamma_h\).}
Introduce the linear functionals
\[
\ell_u(\tauh):=b_{u\sigma}(\tauh,\uu),
\qquad
\ell_\gamma(\tauh):=c_{\gamma\sigma}(\tauh,\gh),
\qquad
\ell_m:=\ell_u+\ell_\gamma.
\]
With respect to the compliance bilinear form \(a_{\sigma\sigma}(\cdot,\cdot)\),
define the Riesz representatives
\(\underline\rho_u,\underline\rho_\gamma,\underline\rho_m\in\underline\Sigma_h\)
by
\begingroup
\begin{equation}\label{eq:egourhogammarhom-def}
	a_{\sigma\sigma}(\underline\rho_u,\tauh)=\ell_u(\tauh),\;
	a_{\sigma\sigma}(\underline\rho_\gamma,\tauh)=\ell_\gamma(\tauh),\;
	a_{\sigma\sigma}(\underline\rho_m,\tauh)=\ell_m(\tauh),
	\; \forall\tauh\in\underline\Sigma_h.
\end{equation}
\endgroup
By linearity and uniqueness of the Riesz representation, \(\underline\rho_m=\underline\rho_u+\underline\rho_\gamma\).

\begingroup
\begin{lemma}
	\label{lem:compliance-representatives-matrix-energy}
	For every \(\mm=(\uc,\gammac)\),
	\[
	(\Km\mm,\mm)=\norm{\underline\rho_m}_{\Ac}^2+s_\gamma(\gh,\gh),
	\qquad
	(\Dm\mm,\mm)=\norm{\underline\rho_u}_{\Ac}^2+\norm{\underline\rho_\gamma}_{\Ac}^2+s_\gamma(\gh,\gh).
	\]
\end{lemma}

\begin{proof}
	Let \({\boldsymbol y_u}\), \({\boldsymbol y_\gamma}\), and \({\boldsymbol y_m}\) be the coefficient vectors of \(\underline\rho_u\), \(\underline\rho_\gamma\), and \(\underline\rho_m\), respectively. By the block definitions \eqref{eq:block-BC-0527}, the coefficient vectors of \(\ell_u\), \(\ell_\gamma\), and \(\ell_m\) are \({\boldsymbol b_u}=\B^{\top}\uc\), \({\boldsymbol b_\gamma}=\C^{\top}\gammac\), and \({\boldsymbol b_m}=\B^{\top}\uc+\C^{\top}\gammac\). Thus the above Riesz relations are equivalent to \(\A {\boldsymbol y_u}={\boldsymbol b_u}\), \(\A {\boldsymbol y_\gamma}={\boldsymbol b_\gamma}\), and \(\A {\boldsymbol y_m}={\boldsymbol b_m}\). It follows that
	\[
	\norm{\underline\rho_m}_{\Ac}^2={\boldsymbol y_m}^{\top}\A {\boldsymbol y_m}={\boldsymbol b_m}^{\top}\A^{-1}{\boldsymbol b_m}=(\B^{\top}\uc+\C^{\top}\gammac)^{\top}\A^{-1}(\B^{\top}\uc+\C^{\top}\gammac)=(\Kmzero\mm,\mm).
	\]
	Since \(\Km=\Kmzero+\diag(0,\Jg)\), adding \(s_\gamma(\gh,\gh)=(\Jg\gammac,\gammac)\) gives the first identity of the lemma. Applying the same argument to \({\boldsymbol b_u}\) and \({\boldsymbol b_\gamma}\) separately gives
	\[
	\norm{\underline\rho_u}_{\Ac}^2+\norm{\underline\rho_\gamma}_{\Ac}^2=(\Kuu\uc,\uc)+(\Kggzero\gammac,\gammac).
	\]
	Adding the same stabilization energy gives the identity for \(\Dm\).
\end{proof}
\endgroup

The next estimate for the rotation representative is the key parameter-robust step. It uses only the deviatoric part of the rotation representative and therefore avoids the degenerating spherical weight in the compliance decomposition above.

\begin{lemma}
	\label{lem:rotation-representative-deviatoric-bound}
	There exists a constant \(C_{\rm as}>0\), independent of \(h\), \(\lambda\), and \(\mu\), such that
	\[
	\norm{\underline\rho_\gamma}_{\Ac}\le C_{\rm as}\sqrt{2\mu}\,\norm{\gh}_{L^2(\Omega)}.
	\]
\end{lemma}

\begin{proof}
	{Taking \(\tauh=\underline\rho_\gamma\) in \eqref{eq:egourhogammarhom-def}} gives \(\norm{\underline\rho_\gamma}_{\Ac}^2=(\asop(\underline\rho_\gamma),\gh)\). Since \(\asop(q\tI)=0\) for every scalar \(q\), the rotation functional annihilates the spherical part, and hence \((\asop(\underline\rho_\gamma),\gh)=(\asop(\underline\rho_\gamma^D),\gh)\). Using the Cauchy--Schwarz inequality,
	we obtain
	$\norm{\underline\rho_\gamma}_{\Ac}^2\le C_{\rm as}\norm{\underline\rho_\gamma^D}_{L^2(\Omega)}\norm{\gh}_{L^2(\Omega)}
	\le C_{\rm as}\sqrt{2\mu}\,\norm{\underline\rho_\gamma}_{\Ac}\norm{\gh}_{L^2(\Omega)}$.
	Dividing by \(\norm{\underline\rho_\gamma}_{\Ac}\) when it is nonzero proves the estimate; the zero case is immediate.
\end{proof}

\begin{theorem}
	\label{thm:main-mech}
	There exist constants \(c_m,C_m>0\), independent of \(h\), \(\lambda\), and \(\mu\), such that
	\[
	c_m(\Dm\mm,\mm)\le(\Km\mm,\mm)\le C_m(\Dm\mm,\mm),\qquad \forall\mm=(\uc,\gammac).
	\]
\end{theorem}
\begin{proof}
	By Lemma~\ref{lem:compliance-representatives-matrix-energy},
	\[
	(\Km\mm,\mm)=\norm{\underline\rho_u+\underline\rho_\gamma}_{\Ac}^2+s_\gamma(\gh,\gh)
	\le2\norm{\underline\rho_u}_{\Ac}^2+2\norm{\underline\rho_\gamma}_{\Ac}^2+s_\gamma(\gh,\gh)
	\le2(\Dm\mm,\mm).
	\]
	For the lower bound, Lemma~\ref{lem:method2-stabilized-elasticity-stability}, the definition of \(\ell_m\), and \(\|{\underline\tau_h}\|_{\Ac}\le(2\mu)^{-1/2}\|{\underline\tau_h}\|_{L^2{(\Omega)}}\) give
	$
	c_E^2\|(\uu,\gh)\|_h^2
	\le
	\frac1{2\mu}\|\underline\rho_m\|_{\Ac}^2
	+|\gh|_{\gamma,h}^2.
	$
	Since \(s_\gamma(\gh,\gh)=\delta_\gamma\mu|\gh|_{\gamma,h}^2\), it follows, for fixed \(\delta_\gamma>0\), that
	$	2\mu\|\gh\|_{L^2(\Omega)}^2
	\lesssim
	\|\underline\rho_m\|_{\Ac}^2+s_\gamma(\gh,\gh)
	=(\Km\mm,\mm)$.
	Lemma~\ref{lem:rotation-representative-deviatoric-bound} then controls \(\|\underline\rho_\gamma\|_{\Ac}^2\) by \((\Km\mm,\mm)\); the identity \(\underline\rho_u=\underline\rho_m-\underline\rho_\gamma\) controls \(\|\underline\rho_u\|_{\Ac}^2\), while \(s_\gamma(\gh,\gh)\le(\Km\mm,\mm)\) holds directly. Hence
	$(\Dm\mm,\mm)
	=\|\underline\rho_u\|_{\Ac}^2+\|\underline\rho_\gamma\|_{\Ac}^2+s_\gamma(\gh,\gh)
	\lesssim(\Km\mm,\mm).
	$
	This gives the required equivalence.
\end{proof}
Finally, we note that the stabilized rotation block is uniformly mass-equivalent.
\begin{corollary}
	\label{cor:Kgg-mass-equivalence}
	There exist \(c_\gamma,C_\gamma>0\), independent of \(h,\lambda,\mu\), such that \(c_\gamma\,\mu\Mg\preceq\Kgg\preceq C_\gamma\,\mu\Mg\).
\end{corollary}
\begin{proof}
	Let \(\gh\in\Gamma_h\) be the macrocell rotation associated with \(\gammac\), and let \(\underline\rho_\gamma\) be its compliance-Riesz representative. Taking \(\mm=(0,\gammac)\) in Lemma~\ref{lem:compliance-representatives-matrix-energy} gives \(E_\gamma:=(\Kgg\gammac,\gammac)=\|\underline\rho_\gamma\|_{\Ac}^2+s_\gamma(\gh,\gh)=\|\underline\rho_\gamma\|_{\Ac}^2+\delta_\gamma\mu|\gh|_{\gamma,h}^2\). For the upper bound, Lemma~\ref{lem:rotation-representative-deviatoric-bound} gives \(\|\underline\rho_\gamma\|_{\Ac}^2\le2C_{\rm as}^2\mu\|\gh\|_{L^2(\Omega)}^2\), while the face inverse estimate gives
	$|\gh|_{\gamma,h}^2=\sum_{F\in\mathcal F_M^{\rm int}}h_F\|\jump{\gh}\|_{L^2(F)}^2\le C_J\|\gh\|_{L^2(\Omega)}^2.$
	Hence \(E_\gamma\le(2C_{\rm as}^2+\delta_\gamma C_J)\mu\|\gh\|_{L^2(\Omega)}^2\).
	For the lower bound, the rotation-energy estimate established in the proof of Theorem~\ref{thm:main-mech}, applied with \(\mm=(0,\gammac)\), gives \(2\mu\|\gh\|_{L^2(\Omega)}^2\lesssim(\Km\mm,\mm)=E_\gamma\).
	Since \((\Mg\gammac,\gammac)=\|\gh\|_{L^2(\Omega)}^2\), the two bounds prove the equivalence.
\end{proof}

\subsection{The uncoupled preconditioners}
\label{subsec:three-preconditioners}
Condensation gives the natural pressure block \(\Tp=\Mp+(\App-\Spsig)+\Spz\), equivalently \(\Tp=c_0\Mzero+(\clambda\Mzero-\Spsig)+\Spz\), whereas the simpler positive definite pressure block is \(\Pmz:=\Mp+\App+\Spz=\Tp+\Spsig\). We consider the two block-diagonal preconditioners
\begin{equation}\label{eq:uncoupled-block-preconditioners}
	\mathcal P_{\rm bd}^{\Tp}=\diag(\Dm,\Tp),
	\qquad
	\mathcal P_{\rm bd}^{\Pmz}=\diag(\Dm,\Pmz).
\end{equation}
We refer to these two preconditioners as uncoupled preconditioners, since their pressure blocks do not involve the mechanics-induced Schur correction.
We first determine when the natural pressure block can be inverted. Define the normalized coercivity constant of the non-storage pressure part
\[
\delta_h
:=
\inf_{\pp\ne0}
\frac{\bigl(((\App-\Spsig)+\Spz)\pp,\pp\bigr)}
{(\Mzero\pp,\pp)}.
\]
Hence, \(\Tp\succ0\) if and only if \(c_0+\delta_h>0\).
For the remainder of this subsection, we assume \(\Tp\succ0\).

To quantify separately the relative size of the stress-induced pressure block and the omitted mechanics--pressure coupling in the \(\diag(\Dm,\Tp)\)-scaled norm, introduce the following two dimensionless quantities:
\begingroup
\small
\begin{equation*}
	\begin{aligned}
		\Thetah
		:=\left\|\Tp^{-1/2}\Spsig\Tp^{-1/2}\right\|_2
		=\sup_{\pp\ne0}
		\frac{(\Spsig\pp,\pp)}{(\Tp\pp,\pp)},
		{\Gammah}
		:=\left\|\Dm^{-1/2}\Lop\Tp^{-1/2}\right\|_2
		=\sup_{\substack{\nn\ne0\\\pp\ne0}}
		\frac{|(\Lop\pp,\nn)|}
		{(\Dm\nn,\nn)^{1/2}(\Tp\pp,\pp)^{1/2}}.
	\end{aligned}
\end{equation*}
\endgroup
Thus \(\Thetah\) is the smallest constant for which
\begin{equation}
	\label{eq:Spsig-controlled-by-Theta-Tp}
	\Spsig\preceq\Thetah\Tp.
\end{equation}
The next lemma connects these two quantities.
\begin{lemma}
	\label{lem:coupling-estimates}
	For every pressure vector \(\pp\) and mechanics vector \(\nn\), 		\begin{equation}
		\label{eq:coupling-estimate-Km}
		\abs{(\Lop\pp,\nn)}
		\le (\Spsig\pp,\pp)^{1/2}(\Km\nn,\nn)^{1/2}.
	\end{equation}
	Consequently,
	\begin{equation}
		\label{eq:coupling-estimate-Tp}
		\abs{(\Lop\pp,\nn)}
		\le \sqrt{C_m\Thetah}\,
		(\Tp\pp,\pp)^{1/2}(\Dm\nn,\nn)^{1/2},
		\qquad
		{\Gammah^2}\le C_m\Thetah.
	\end{equation}
\end{lemma}

\begin{proof}
	{Writing \(\nn=(\mathbf v,\boldsymbol\eta)\) and using the definition of \(\Lop\), we obtain}
 	$(\Lop\pp,\nn)
	=(\A^{-1/2}\E\pp,
	\A^{-1/2}(\B^{\top}\mathbf v+\C^{\top}\boldsymbol\eta))$.
	Cauchy--Schwarz gives
	$|(\Lop\pp,\nn)|
	\le(\Spsig\pp,\pp)^{1/2}(\Kmzero\nn,\nn)^{1/2}
	\le(\Spsig\pp,\pp)^{1/2}(\Km\nn,\nn)^{1/2},$
	because \(\Km=\Kmzero+\diag(0,\Jg)\succeq\Kmzero\). This proves \eqref{eq:coupling-estimate-Km}. Combining it with \(\Km\preceq C_m\Dm\) and \eqref{eq:Spsig-controlled-by-Theta-Tp} gives \eqref{eq:coupling-estimate-Tp}.
\end{proof}

Since \(\Pmz=\Tp+\Spsig\), relation \eqref{eq:Spsig-controlled-by-Theta-Tp} yields $\Tp\preceq\Pmz\preceq(1+\Thetah)\Tp$.
Introduce the energy norms induced by the two preconditioners,
\[
\norm{(\mm,\pp)}_{\DT}^2
:=(\Dm\mm,\mm)+(\Tp\pp,\pp),
\qquad
\norm{(\mm,\pp)}_{\DPmz}^2
:=(\Dm\mm,\mm)+(\Pmz\pp,\pp).
\]
Then
\begin{equation}
	\label{eq:DT-DPmz-norm-equivalence}
	\norm{(\mm,\pp)}_{\DT}^2
	\le \norm{(\mm,\pp)}_{\DPmz}^2
	\le(1+\Thetah)\norm{(\mm,\pp)}_{\DT}^2.
\end{equation}
The bilinear form associated with \(\Ared\) is
\[
\mathcal B_{\rm red}((\mm,\pp),(\nn,\qq))
:=
-(\Km\mm,\nn)
-(\Lop\pp,\nn)
-(\Lop\qq,\mm)
+(\Tp\pp,\qq).
\]
Finally, set
$C_T(\theta)
:=\frac{C_m+1+\sqrt{(C_m-1)^2+4C_m\theta}}{2}$,
and
$\beta_T:=\min\{c_m,1\}$.

\begin{theorem}
	\label{thm:uncoupled-Theta-main}
	For every parameter choice satisfying \(\Tp\succ0\), the preconditioner \(\mathcal P_{\rm bd}^{\Tp}\) satisfies
	\begin{align}
		\abs{\mathcal B_{\rm red}((\mm,\pp),(\nn,\qq))}
		&\le C_T(\Thetah)
		\norm{(\mm,\pp)}_{\DT}\norm{(\nn,\qq)}_{\DT},
		\label{eq:Tp-continuity-main}\\
		\inf_{(\mm,\pp)\ne0}\sup_{(\nn,\qq)\ne0}
		\frac{\mathcal B_{\rm red}((\mm,\pp),(\nn,\qq))}
		{\norm{(\mm,\pp)}_{\DT}\norm{(\nn,\qq)}_{\DT}}
		&\ge\beta_T.
		\label{eq:Tp-infsup-main}
	\end{align}
	The preconditioner \(\mathcal P_{\rm bd}^{\Pmz}\) satisfies
	\begin{align*}
		\abs{\mathcal B_{\rm red}((\mm,\pp),(\nn,\qq))}
		&\le C_T(\Thetah)
		\norm{(\mm,\pp)}_{\DPmz}\norm{(\nn,\qq)}_{\DPmz},
		\\
		\inf_{(\mm,\pp)\ne0}\sup_{(\nn,\qq)\ne0}
		\frac{\mathcal B_{\rm red}((\mm,\pp),(\nn,\qq))}
		{\norm{(\mm,\pp)}_{\DPmz}\norm{(\nn,\qq)}_{\DPmz}}
		&\ge\frac{\beta_T}{1+\Thetah}.
	\end{align*}
\end{theorem}

\begin{proof}
	Set \(a=(\Dm\mm,\mm)^{1/2}\), \(b=(\Tp\pp,\pp)^{1/2}\), \(c=(\Dm\nn,\nn)^{1/2}\),\newline and \(d=(\Tp\qq,\qq)^{1/2}\). Theorem~\ref{thm:main-mech} and Lemma~\ref{lem:coupling-estimates} yield
	\[
	\abs{\mathcal B_{\rm red}((\mm,\pp),(\nn,\qq))}
	\le
	\begin{bmatrix}a&b\end{bmatrix}
	\begin{bmatrix}
		C_m&\sqrt{C_m\Thetah}\\
		\sqrt{C_m\Thetah}&1
	\end{bmatrix}
	\begin{bmatrix}c\\d\end{bmatrix}.
	\]
	The spectral norm of the middle matrix is \(C_T(\Thetah)\), proving \eqref{eq:Tp-continuity-main}. For the inf--sup bound choose \((\nn,\qq)=(-\mm,\pp)\); then
	$\mathcal B_{\rm red}((\mm,\pp),(-\mm,\pp))
	=(\Km\mm,\mm)+(\Tp\pp,\pp)
	\ge\beta_T\norm{(\mm,\pp)}_{\DT}^2$.
	This proves \eqref{eq:Tp-infsup-main}. The two estimates for \(\mathcal P_{\rm bd}^{\Pmz}\) follow from the first two estimates and \eqref{eq:DT-DPmz-norm-equivalence}.
\end{proof}

{The bounds in Theorem~\ref{thm:uncoupled-Theta-main} are expressed in terms of the discrete ratio \(\Thetah\).} We next derive a transparent sufficient bound in terms of the physical parameters. By the definition of \(\delta_h\), $\Tp\succeq(c_0+\delta_h)\Mzero$,
while \(\Spsig\preceq\App=\clambda\Mzero\). Therefore
\begin{equation*}
	\Thetah
	\le
	\frac{\clambda}{c_0+\delta_h}
	=
	\frac{\alpha_B^2}{\Kdr(c_0+\delta_h)}.
\end{equation*}
For \(c_0>0\), define
\begin{equation*}
	\chi_{\lambda,c_0}
	:=
	\frac{\clambda}{c_0}
	=
	\frac{\alpha_B^2}{c_0\Kdr}
	=
	\frac{d\alpha_B^2}{c_0(2\mu+d\lambda)}.
\end{equation*}
Then \(\Thetah\le\chi_{\lambda,c_0}\), and the following parameter-explicit consequence is immediate.
\begin{corollary}
	\label{cor:uncoupled-parameter-explicit}
	Let \(\mathfrak P\) be a family of parameter choices with \(c_0>0\) such that
	\begin{equation}
		\label{eq:balanced-storage-explicit}
		\chi_*:=\sup_{\mathfrak P}
		\frac{\alpha_B^2}{c_0\Kdr}<\infty.
	\end{equation}
	Then \(\mathcal P_{\rm bd}^{\Tp}\) and \(\mathcal P_{\rm bd}^{\Pmz}\) are uniformly stable on \(\mathfrak P\). Their continuity and inf--sup constants depend only on \(c_m\), \(C_m\), and \(\chi_*\), and are independent of \(h\), \(\Delta t\), \(K\), \(\lambda\), and \(c_0\).
\end{corollary}
Condition \eqref{eq:balanced-storage-explicit} is sufficient rather than necessary. In particular, it permits the critical scaling
\[
c_0=C\frac{\alpha_B^2}{\Kdr},
\qquad C>0,
\]
and hence{, for fixed \(\alpha_B>0\) and \(\mu\),} \(c_0\asymp\Kdr^{-1}\sim\lambda^{-1}\) in the nearly incompressible limit. {More generally, Theorem~\ref{thm:uncoupled-Theta-main} yields uniform stability bounds provided that \(\Tp\succ0\) and \(\Thetah\) is uniformly bounded.} Thus, when \(c_0=0\), the uncoupled analysis may still apply if \(\delta_h>0\) and the residual compliance and Darcy terms provide a uniform bound for \(\Thetah\). If \(c_0=\delta_h=0\), then \(\Tp\) is singular: \(\mathcal P_{\rm bd}^{\Tp}\) is not defined on the full pressure space, and the present \(\Thetah\)-based argument does not yield a stability estimate for \(\mathcal P_{\rm bd}^{\Pmz}\).

It is useful to close this subsection by comparing the two uncoupled pressure blocks with the exact pressure Schur complement. Define the mechanics-induced pressure correction \(\mathsf Q:=\Lop^{\top}\Km^{-1}\Lop\) and the exact Schur complement \(\Sp:=\Tp+\mathsf Q\). The following result compares \(\Tp\), \(\Sp\), and \(\Pmz\).
\begin{proposition}
	\label{prop:Q-bounded-by-Spsig}
	The mechanics-induced correction satisfies
	\begin{equation}
		\label{eq:Q-bounded-by-Spsig}
		0\preceq\mathsf Q\preceq\Spsig.
	\end{equation}
	Consequently,
	\begin{equation}
		\label{eq:Tp-Sp-Pmz-ordering}
		\Tp\preceq\Sp\preceq\Pmz.
	\end{equation}
\end{proposition}

\begin{proof}
	Since \(\Km\succ0\), \(\mathsf Q\succeq0\). For every \(\pp\), \eqref{eq:coupling-estimate-Km} gives
	$(\mathsf Q\pp,\pp)
	=\sup_{\nn\ne\mathbf0}
	\frac{|(\Lop\pp,\nn)|^2}{(\Km\nn,\nn)}
	\le(\Spsig\pp,\pp)$.
	Thus \eqref{eq:Q-bounded-by-Spsig} holds; using \(\Sp=\Tp+\mathsf Q\) and \(\Pmz=\Tp+\Spsig\) gives \eqref{eq:Tp-Sp-Pmz-ordering}.
\end{proof}

The relative size of the exact mechanics-induced pressure correction is measured by \(\Theta_h^{\rm sch}:=\left\|\Tp^{-1/2}\mathsf Q\Tp^{-1/2}\right\|_2\). Combining Theorem~\ref{thm:main-mech} with Proposition~\ref{prop:Q-bounded-by-Spsig} gives the hierarchy
\begin{equation}
	\label{eq:three-level-coupling-hierarchy}
	c_m\Theta_h^{\rm sch}
	\le {\Gammah^2}
	\le C_m\Theta_h^{\rm sch}
	\le C_m\Thetah
	\le
	\frac{C_m\alpha_B^2}{\Kdr(c_0+\delta_h)}.
\end{equation}
Indeed, \(c_m\Dm\preceq\Km\preceq C_m\Dm\) implies \((1/C_m)\Dm^{-1}\preceq\Km^{-1}\preceq(1/c_m)\Dm^{-1}\), while \(\mathsf Q\preceq\Spsig\). Thus the three pressure blocks have a nested interpretation: the natural block \(\Tp\) omits the mechanics correction, the exact Schur block \(\Sp\) adds the true correction \(\mathsf Q\), and the simpler block \(\Pmz\) replaces \(\mathsf Q\) by the larger stress-compliance envelope \(\Spsig\). Formula \eqref{eq:three-level-coupling-hierarchy} also separates three levels of information: \(\Theta_h^{\rm sch}\) is the exact Schur-correction ratio, \(\Thetah\) is its stress-compliance envelope, and \(\chi_{\lambda,c_0}\) is the coarser storage-only envelope when \(c_0>0\). Hence growth of \(\chi_{\lambda,c_0}\) alone cannot prove deterioration of the actual coupling; it only shows that storage has ceased to provide a uniform proof. This motivates the Schur-consistent pressure blocks studied next, which approximate the correction \(\mathsf Q\) instead of relying only on the envelope \(\Spsig\).

\subsection{Spectral analysis of the Schur-consistent preconditioners}
\label{subsec:schur-spectral-analysis}

The exact pressure Schur complement is used only as an algebraic benchmark. We first analyze the exact block factors and then replace \(\Km^{-1}\) by the computable block-diagonal inverse \(\Dm^{-1}\).
Recall the exact mechanics-induced pressure correction $\mathsf Q:=\Lop^{\top}\Km^{-1}\Lop$ and $\Sp:=\Tp+\mathsf Q$, and define its Schur-consistent counterpart by
\begin{equation}
	\label{eq:Psc-def}
	\mathsf Q_D:=\Lop^{\top}\Dm^{-1}\Lop,
	\qquad
	\Psc:=\Tp+\mathsf Q_D.
\end{equation}
Since \(\Km,\Dm\succ0\), \(\Tp\succeq0\),
	and the condensed matrix \(\Ared\) is nonsingular by
	Theorem~\ref{thm:fully-discrete-wellposedness},
	both \(\Sp\) and \(\Psc\) are positive definite.
Define the exact and Schur-consistent block-diagonal preconditioners by
\begin{equation}\label{eq:schur-block}
	\mathcal P_{\rm bd}^{\rm ex}
	:=\begin{bmatrix}\Km&0\\0&\Sp\end{bmatrix},
	\qquad
	\mathcal P_{\rm bd}^{\rm sc}
	:=\begin{bmatrix}\Dm&0\\0&\Psc\end{bmatrix},
\end{equation}
and the lower-triangular factors by
\begin{equation}\label{eq:schur-tri-block}
	\mathcal P_{\rm lt}^{\rm ex}
	:=\begin{bmatrix}\Km&0\\\Lop^{\top}&-\Sp\end{bmatrix},
	\qquad
	\mathcal P_{\rm lt}^{\rm sc}
	:=\begin{bmatrix}\Dm&0\\\Lop^{\top}&-\Psc\end{bmatrix}.
\end{equation}
For the sign convention of \(\Ared\), set
$		\mathcal P_{{\rm lt},{\rm red}}^{\rm ex}:=-\mathcal P_{\rm lt}^{\rm ex}$,
	and	$\mathcal P_{{\rm lt},{\rm red}}^{\rm sc}:=-\mathcal P_{\rm lt}^{\rm sc}$.

The following result is the bridge from the exact Schur complement to its computable counterpart. It shows that replacing $\Km^{-1}$ by $\Dm^{-1}$ in the mechanics-induced pressure correction preserves the pressure energy uniformly.
\begin{theorem}
	\label{thm:Psc-equivalent-Sp}
	Let \(c_m\) and \(C_m\) be the constants in Theorem~\ref{thm:main-mech}, and set $\alpha_{\rm sc}:=\min\{1,C_m^{-1}\}$ and $\beta_{\rm sc}:=\max\{1,c_m^{-1}\}$. Then
	\begin{equation}
		\label{eq:Psc-Sp-equivalence}
		\alpha_{\rm sc}\Psc\preceq\Sp\preceq\beta_{\rm sc}\Psc.
	\end{equation}
	\end{theorem}

\begin{proof}
	From~Theorem~\ref{thm:main-mech}, we have $C_m^{-1}\Dm^{-1}\preceq\Km^{-1}\preceq c_m^{-1}\Dm^{-1}$. Addition of \(\Tp\succeq0\) yields $\Tp+C_m^{-1}\mathsf Q_D \preceq\Sp \preceq \Tp+c_m^{-1}\mathsf Q_D$, which proves \eqref{eq:Psc-Sp-equivalence} since \(\Tp\succeq0\) and \(\mathsf Q_D\succeq0\).
\end{proof}
We next establish eigenvalue bounds for the ideal block-diagonal preconditioner of the reduced system, which contains the positive semidefinite pressure block $\Tp$.
\begin{theorem}
	\label{thm:exact-schur-for-Ared}
	Let \(\varphi=(1+\sqrt5)/2\). Then
	\begin{equation}
		\label{eq:Ared-exact-bd-spectrum}
		\spec\bigl((\mathcal P_{\rm bd}^{\rm ex})^{-1}\Ared\bigr)
		\subset[-\varphi,-1]\cup[\varphi^{-1},1].
	\end{equation}
	{Consequently}, for every vector \(\boldsymbol x\),
	\begin{equation}
		\label{eq:Ared-exact-bd-dual-norm-equivalence}
		\varphi^{-1}\|\boldsymbol x\|_{\mathcal P_{\rm bd}^{\rm ex}}
		\le
		\|\Ared\boldsymbol x\|_{(\mathcal P_{\rm bd}^{\rm ex})^{-1}}
		\le
		\varphi\|\boldsymbol x\|_{\mathcal P_{\rm bd}^{\rm ex}}.
	\end{equation}
\end{theorem}

\begin{proof}
	Set $\mathsf H:=\Km^{-1/2}\Lop\Sp^{-1/2}$, so that
	\begin{equation}
		\label{eq:H-contraction}
		0\preceq \mathsf H^{\top}\mathsf H
		=\Sp^{-1/2}\Lop^{\top}\Km^{-1}\Lop\Sp^{-1/2}
		\preceq I.
	\end{equation}
	The preconditioned matrix is similar to the symmetric matrix
	\begin{equation*}
		{\mathsf T_{\rm bd}^{\rm ex}}
		:=(\mathcal P_{\rm bd}^{\rm ex})^{-1/2}
		\Ared
		(\mathcal P_{\rm bd}^{\rm ex})^{-1/2}
		=
		\begin{bmatrix}
			-I&-\mathsf H\\
			-\mathsf H^{\top}&I-\mathsf H^{\top}\mathsf H
		\end{bmatrix}.
	\end{equation*}
	Let \(\sigma_i\in(0,1]\) be a positive singular value of \(\mathsf H\), with corresponding left and right singular vectors \(\boldsymbol q_i\) and \(\boldsymbol p_i\), and put \(\omega_i:=\sigma_i^2\). The subspace spanned by \((\boldsymbol q_i,0)^{\top}\) and \((0,\boldsymbol p_i)^{\top}\) is invariant under \({\mathsf T_{\rm bd}^{\rm ex}}\), and the restriction of \({\mathsf T_{\rm bd}^{\rm ex}}\) to this subspace is
	\[
	\begin{bmatrix}
		-1&-\sigma_i\\
		-\sigma_i&1-\sigma_i^2
	\end{bmatrix}.
	\]
	Its eigenvalues are $\rho_i^{\pm}=\frac{-\omega_i\pm\sqrt{\omega_i^2+4}}{2}$. For \(0<\omega_i\le1\), the negative branch lies in \([-\varphi,-1)\), and the positive branch lies in {\([\varphi^{-1},1)\)}. It remains to treat the case \(\sigma_i=0\). For \(\mathsf H^{\top}\boldsymbol q=0\), we have \({\mathsf T_{\rm bd}^{\rm ex}}(\boldsymbol q,0)^{\top}=-(\boldsymbol q,0)^{\top}\), and if \(\mathsf H\boldsymbol p=0\), then \({\mathsf T_{\rm bd}^{\rm ex}}(0,\boldsymbol p)^{\top}=(0,\boldsymbol p)^{\top}\). Therefore, the remaining directions associated with zero singular values contribute only the eigenvalues $1$ and $-1$. This proves \eqref{eq:Ared-exact-bd-spectrum}.

	Because \({\mathsf T_{\rm bd}^{\rm ex}}\) is symmetric, \eqref{eq:Ared-exact-bd-spectrum} implies $\varphi^{-1}\|\boldsymbol y\|_2 \le\|{\mathsf T_{\rm bd}^{\rm ex}}\boldsymbol y\|_2 \le\varphi\|\boldsymbol y\|_2$. Taking \(\boldsymbol y=(\mathcal P_{\rm bd}^{\rm ex})^{1/2}\boldsymbol x\) gives \eqref{eq:Ared-exact-bd-dual-norm-equivalence}.
\end{proof}

Combining the exact benchmark with the uniform equivalences $\Km\simeq\Dm$ and $\Sp\simeq\Psc$ transfers the same parameter-independent spectral separation to the Schur-consistent block-diagonal preconditioner.
\begin{corollary}
	\label{cor:Pbd-sc-ex-equivalence}
	Set $a_{\rm bd}:=\min\{c_m,\alpha_{\rm sc}\}$, $b_{\rm bd}:=\max\{C_m,\beta_{\rm sc}\}$. Then
	\begin{equation}
		\label{eq:Ared-sc-bd-spectrum}
		\spec\bigl((\mathcal P_{\rm bd}^{\rm sc})^{-1}\Ared\bigr)
		\subset
		[-\varphi b_{\rm bd},-a_{\rm bd}/\varphi]
		\cup
		[a_{\rm bd}/\varphi,\varphi b_{\rm bd}].
	\end{equation}
\end{corollary}

\begin{proof}
	The block inequalities in Theorem~\ref{thm:main-mech} and Theorem~\ref{thm:Psc-equivalent-Sp} give \(a_{\rm bd}\mathcal P_{\rm bd}^{\rm sc}\preceq\mathcal P_{\rm bd}^{\rm ex}\preceq b_{\rm bd}\mathcal P_{\rm bd}^{\rm sc}\). Hence \(\sqrt{a_{\rm bd}}\|\boldsymbol x\|_{\mathcal P_{\rm bd}^{\rm sc}}\le\|\boldsymbol x\|_{\mathcal P_{\rm bd}^{\rm ex}}\le\sqrt{b_{\rm bd}}\|\boldsymbol x\|_{\mathcal P_{\rm bd}^{\rm sc}}\), whereas inversion of the same block equivalence gives \(\sqrt{a_{\rm bd}}\|\boldsymbol v\|_{(\mathcal P_{\rm bd}^{\rm ex})^{-1}}\le\|\boldsymbol v\|_{(\mathcal P_{\rm bd}^{\rm sc})^{-1}}\le\sqrt{b_{\rm bd}}\|\boldsymbol v\|_{(\mathcal P_{\rm bd}^{\rm ex})^{-1}}\). Combining these estimates with \eqref{eq:Ared-exact-bd-dual-norm-equivalence} yields \(\frac{a_{\rm bd}}{\varphi}\|\boldsymbol x\|_{\mathcal P_{\rm bd}^{\rm sc}}\le\|\Ared\boldsymbol x\|_{(\mathcal P_{\rm bd}^{\rm sc})^{-1}}\le\varphi b_{\rm bd}\|\boldsymbol x\|_{\mathcal P_{\rm bd}^{\rm sc}}\). The matrix \((\mathcal P_{\rm bd}^{\rm sc})^{-1/2}\Ared(\mathcal P_{\rm bd}^{\rm sc})^{-1/2}\) is symmetric, so the same norm bounds imply \eqref{eq:Ared-sc-bd-spectrum}.
\end{proof}

A direct block multiplication gives the exact factorization identity
\begin{equation}
	\label{eq:Ared-exact-lt-preconditioned}
	(\mathcal P_{{\rm lt},{\rm red}}^{\rm ex})^{-1}\Ared
	=
	\begin{bmatrix}
		I&\Km^{-1}\Lop\\
		0&I
	\end{bmatrix}.
\end{equation}
In particular, every eigenvalue is equal to one and $\bigl((\mathcal P_{{\rm lt},{\rm red}}^{\rm ex})^{-1}\Ared-I\bigr)^2=0$.

{With exact block solves, \eqref{eq:Ared-exact-lt-preconditioned} implies that GMRES converges in at most two iterations in exact arithmetic. To analyze the exact and Schur-consistent triangular factors in a common framework, we derive weighted FOV and operator-norm bounds.} For any SPD matrix \(\mathcal N\), write
$	(\boldsymbol x,\boldsymbol y)_{\mathcal N}
:=\boldsymbol x^{\top}\mathcal N\boldsymbol y$,
and
$\|\boldsymbol x\|_{\mathcal N}^2
:=(\boldsymbol x,\boldsymbol x)_{\mathcal N}$.
\begin{theorem}
	\label{thm:exact-lt-fov}
	Let \(\mathcal N_{\rm ex}:=\mathcal P_{\rm bd}^{\rm ex}\) and \(\mathcal L_{\rm ex}:=(\mathcal P_{{\rm lt},{\rm red}}^{\rm ex})^{-1}\). Then, for every nonzero vector \(\boldsymbol x\),
	\[
	\frac{(\mathcal L_{\rm ex}\Ared\boldsymbol x,\boldsymbol x)_{\mathcal N_{\rm ex}}}
	{\|\boldsymbol x\|_{\mathcal N_{\rm ex}}^2}
	\ge\frac12,
	\qquad
	\frac{\|\mathcal L_{\rm ex}\Ared\boldsymbol x\|_{\mathcal N_{\rm ex}}}
	{\|\boldsymbol x\|_{\mathcal N_{\rm ex}}}
	\le\varphi.
	\]
\end{theorem}

\begin{proof}
	With \(\mathsf H\) from \eqref{eq:H-contraction} and \(\boldsymbol y=\mathcal N_{\rm ex}^{1/2}\boldsymbol x=(\boldsymbol y_m,\boldsymbol y_p)^{\top}\), set \({\mathsf T_{\rm lt}^{\rm ex}}:=\mathcal N_{\rm ex}^{1/2}\mathcal L_{\rm ex}\Ared\mathcal N_{\rm ex}^{-1/2}\). By \eqref{eq:Ared-exact-lt-preconditioned}, \({\mathsf T_{\rm lt}^{\rm ex}}=\bigl[\begin{smallmatrix}I&\mathsf H\\0&I\end{smallmatrix}\bigr]\) and \(\|\mathsf H\|_2\le1\). Set \(a=\|\boldsymbol y_m\|_2\) and \(b=\|\boldsymbol y_p\|_2\). Then
	\[
	\begin{aligned}
		(\mathcal L_{\rm ex}\Ared\boldsymbol x,\boldsymbol x)_{\mathcal N_{\rm ex}}
		&=a^2+(\mathsf H\boldsymbol y_p,\boldsymbol y_m)+b^2
		\ge a^2-ab+b^2
		\ge\tfrac12(a^2+b^2)
		=\tfrac12\|\boldsymbol x\|_{\mathcal N_{\rm ex}}^2,\\
		\|\mathcal L_{\rm ex}\Ared\boldsymbol x\|_{\mathcal N_{\rm ex}}^2
		&=\|\boldsymbol y_m+\mathsf H\boldsymbol y_p\|_2^2+b^2
		\le(a+b)^2+b^2
		\le\varphi^2(a^2+b^2)
		=\varphi^2\|\boldsymbol x\|_{\mathcal N_{\rm ex}}^2.
	\end{aligned}
	\]
	The second estimate follows by taking square roots.
\end{proof}

Before estimating the Schur-consistent triangularly preconditioned operator itself, we record that the exact and Schur-consistent lower-triangular factors are uniformly comparable at the spectral level.

\begin{corollary}
	\label{cor:Plt-sc-ex-equivalence}
	The exact and Schur-consistent lower-triangular factors satisfy
	\begin{equation*}
		\spec\bigl((\mathcal P_{\rm lt}^{\rm sc})^{-1}\mathcal P_{\rm lt}^{\rm ex}\bigr)
		\subset
		[\min\{c_m,\alpha_{\rm sc}\},\max\{C_m,\beta_{\rm sc}\}].
	\end{equation*}
\end{corollary}

\begin{proof}
	The product \((\mathcal P_{\rm lt}^{\rm sc})^{-1}\mathcal P_{\rm lt}^{\rm ex}\) is block lower triangular with diagonal blocks \(\Dm^{-1}\Km\) and \(\Psc^{-1}\Sp\), so its spectrum is the union of their spectra. Theorem~\ref{thm:main-mech} and Theorem~\ref{thm:Psc-equivalent-Sp} bound these spectra by \([c_m,C_m]\) and \([\alpha_{\rm sc},\beta_{\rm sc}]\), respectively, proving the stated inclusion.
	\end{proof}

The relative spectral comparison of the two triangular factors does not by itself yield a GMRES estimate; the following theorem therefore derives a direct FOV bound for the Schur-consistent factor. Let \(c_m\) be as in Theorem~\ref{thm:main-mech}, and set \(\mathcal N_{\rm sc}:=\mathcal P_{\rm bd}^{\rm sc}\) and \(\mathcal L_{\rm sc}:=(\mathcal P_{{\rm lt},{\rm red}}^{\rm sc})^{-1}\).
\begin{theorem}\relax
	\label{thm:sc-lt-uniform-fov}
	For every nonzero vector \(\boldsymbol x\),
	\begin{equation}
		\label{eq:sc-lt-uniform-fov}
		\frac{(\mathcal L_{\rm sc}\Ared\boldsymbol x,\boldsymbol x)_{\mathcal N_{\rm sc}}}
		{\|\boldsymbol x\|_{\mathcal N_{\rm sc}}^2}
		\ge\Sigma_{\rm sc},
		\qquad
		\frac{\|\mathcal L_{\rm sc}\Ared\boldsymbol x\|_{\mathcal N_{\rm sc}}}
		{\|\boldsymbol x\|_{\mathcal N_{\rm sc}}}
		\le3,
	\end{equation}
	where \(\Sigma_{\rm sc}:=(1-1/\sqrt2)\min\{c_m,1\}>0\).
\end{theorem}

\begin{proof}
	Set \(\mathsf R:=\Dm^{-1/2}\Km\Dm^{-1/2}\) and \(\mathsf G:=\Dm^{-1/2}\Lop\Psc^{-1/2}\). The proof of Theorem~\ref{thm:main-mech} gives the sharper upper bound \(\Km\preceq2\Dm\), so
	$c_m I\preceq \mathsf R\preceq2I$.
	Moreover, \(\Tp\succeq0\) and \eqref{eq:Psc-def} imply
	$\mathsf G^{\top}\mathsf G=\Psc^{-1/2}\mathsf Q_D\Psc^{-1/2}\preceq I$,
	and
	$\mathsf G^{\top}\mathsf R\mathsf G\preceq2\mathsf G^{\top}\mathsf G\preceq2I$.
	A block calculation yields
	${\mathsf T_{\rm lt}^{\rm sc}}:=\mathcal N_{\rm sc}^{1/2}\mathcal L_{\rm sc}\Ared\mathcal N_{\rm sc}^{-1/2}
	=\begin{bmatrix}\mathsf R&\mathsf G\\\mathsf G^{\top}(\mathsf R-I)&I\end{bmatrix}$.

	For \(\boldsymbol y=\mathcal N_{\rm sc}^{1/2}\boldsymbol x=(\boldsymbol y_m,\boldsymbol y_p)^{\top}\), let \(a^2=(\mathsf R\boldsymbol y_m,\boldsymbol y_m)\) and \(b=\|\boldsymbol y_p\|_2\). Since \(\mathsf R\) is symmetric, we have
	$(\mathsf G\boldsymbol y_p,\boldsymbol y_m)
	+(\mathsf G^{\top}(\mathsf R-I)\boldsymbol y_m,\boldsymbol y_p)
	=(\mathsf G\boldsymbol y_p,\mathsf R\boldsymbol y_m)$.
	The Cauchy--Schwarz inequality in the \(\mathsf R\)-inner product therefore gives
	\[
	\bigl|(\mathsf G\boldsymbol y_p,\mathsf R\boldsymbol y_m)\bigr|
	\le a\,(\mathsf G^{\top}\mathsf R\mathsf G\boldsymbol y_p,\boldsymbol y_p)^{1/2}
	\le\sqrt2\,ab.
	\]
	Consequently, we obtain
	$(\mathcal L_{\rm sc}\Ared\boldsymbol x,\boldsymbol x)_{\mathcal N_{\rm sc}}
	\ge\left(1-\frac1{\sqrt2}\right)(a^2+b^2)
	\ge\Sigma_{\rm sc}\|\boldsymbol x\|_{\mathcal N_{\rm sc}}^2$,
	where the last inequality follows from \(\mathsf R\succeq c_m I\) and \(\|\boldsymbol y\|_2=\|\boldsymbol x\|_{\mathcal N_{\rm sc}}\).

	For the operator norm, \(0\prec \mathsf R\preceq2I\) gives \(\|\mathsf R-I\|_2\le1\), while \(\mathsf G^{\top}\mathsf G\preceq I\) gives \(\|\mathsf G\|_2\le1\). Therefore, we obtain
	\begingroup
	\small
	\[
	\begin{aligned}
		\|{\mathsf T_{\rm lt}^{\rm sc}}\|_2&
		\le\left\|\begin{bmatrix}\mathsf R&0\\0&I\end{bmatrix}\right\|_2
		+\left\|\begin{bmatrix}0&\mathsf G\\\mathsf G^{\top}(\mathsf R-I)&0\end{bmatrix}\right\|_2\\
		&=\max\{\|\mathsf R\|_2,1\}
		+\max\{\|\mathsf G\|_2,\|\mathsf G^{\top}(\mathsf R-I)\|_2\}
		\le3.
	\end{aligned}
	\]
	\endgroup
	$\|\mathcal L_{\rm sc}\Ared\boldsymbol x\|_{\mathcal N_{\rm sc}}
	=\|{\mathsf T_{\rm lt}^{\rm sc}}\boldsymbol y\|_2
	\le3\|\boldsymbol y\|_2
	=3\|\boldsymbol x\|_{\mathcal N_{\rm sc}}$,
	which completes the proof.
\end{proof}

\subsection{Practical implementation of the block preconditioners}
The exact block inverses appearing in the preceding preconditioners are not formed explicitly.  For the stabilized formulation, we use a vertex-patch two-grid approximation for the displacement block, a macrocell mass inverse for the rotation block, and a scalar two-grid approximation for the pressure block.

For the uncoupled preconditioners in \eqref{eq:uncoupled-block-preconditioners}, we use
\[
{\mathsf B_u^{-1}}\approx\Kuu^{-1},
\qquad
{\mathsf B_\gamma^{-1}}\approx(\mu\Mg)^{-1}\approx\Kgg^{-1},
\qquad
{\mathsf B_p^{-1}}\approx\Tp^{-1}
\quad\text{or}\quad
{\mathsf B_p^{-1}}\approx(\Pmz)^{-1}.
\]
{The vector \(\mathcal P^1\) coarse basis comprises componentwise samples of the physical-coordinate modes \(1,x_1,\ldots,x_d\) at macrocell centers within each nonoverlapping vertex-centered aggregate. Each vertex-patch Gauss--Seidel block contains displacement degrees of freedom of macrocells sharing a vertex.} By Corollary~\ref{cor:Kgg-mass-equivalence}, the stabilized rotation block is uniformly mass-equivalent.
For the pressure blocks, we use two-grid solvers with forward/backward
Gauss--Seidel pre-/post-smoothing.

The exact Schur complement \(\Sp\) is used only in the analysis. In computation, the stabilized mechanics inverse is approximated by \(\mathcal B_m^{-1}:=\diag({\mathsf B_u^{-1}},{\mathsf B_\gamma^{-1}})\), and the practical Schur-consistent pressure block is \(\widehat{\Psc}:=\Tp+\Lop^{\top}\mathcal B_m^{-1}\Lop\). After this pressure matrix is assembled, its inverse action is approximated by the scalar pressure two-grid solver.
The lower-triangular preconditioner uses the same block solvers together with one lower-left coupling update in the usual block forward substitution.
\section{Numerical experiments}\label{sec:numerics}
This section presents 2D and 3D numerical tests demonstrating the robustness of the proposed algorithms with respect to physical and discretization parameters. Throughout, we use isotropic permeability \(K=\kappa I\), backward Euler time stepping, and the rotation-jump stabilization parameter \(\delta_\gamma=0.5\). The condensed control-volume three-field systems are solved by preconditioned GMRES with relative residual tolerances \(10^{-8}\) in 2D and \(10^{-6}\) in 3D.
In the iteration tables, an entry \(a\,(b)\) gives the GMRES iteration counts for \(\mathcal P_{\rm bd}^{\Pmz}\) and \(\mathcal P_{\rm bd}^{\rm sc}\), respectively.
We mainly present results for \(\mathcal P_{\mathrm{bd}}^{\Pmz}\) and \(\mathcal P_{\mathrm{bd}}^{\mathrm{sc}}\). For the tested parameter ranges, \(\mathcal P_{\mathrm{bd}}^{\Tp}\) and \(\mathcal P_{\mathrm{lt}}^{\mathrm{sc}}\) yield iteration counts very close to those of \(\mathcal P_{\mathrm{bd}}^{\Pmz}\) and \(\mathcal P_{\mathrm{bd}}^{\mathrm{sc}}\), respectively; hence, the results for the former two are omitted for brevity.

    \begin{example}[2D manufactured solution on rectangular meshes]\label{example:accuracy-manufacture}	\normalfont
		We consider a manufactured solution on \(\Omega=[0,1]^2\). The exact pressure is \(p(x,y,t)=e^{-t}\sin(\pi x)\sin(\pi y)\). The displacement components are \(u_1=\allowbreak\frac{e^{-t}}2[\sin^2(\pi x)\sin(2\pi y)\allowbreak+\lambda^{-1}\sin(2\pi x)\sin^2(\pi y)]\) and \(u_2=\allowbreak\frac{e^{-t}}2[-\sin(2\pi x)\sin^2(\pi y)\allowbreak+\lambda^{-1}\sin^2(\pi x)\sin(2\pi y)]\). The remaining variables \(\gamma\), \(\bs z\), \(\underline\sigma\), and the source terms are obtained from the Biot equations. The final time is \(T=0.04\). The reduced linear systems in the accuracy study are solved with the \(\mathcal P_{\rm bd}^{\Pmz}\) preconditioner.
	Figure~\ref{fig:example1-accuracy-convergence} checks the approximation accuracy and the locking-free behavior of the discretization. All plotted errors are absolute \(L^2\) errors. The macroelement errors \(e_{\mathcal M\sigma}\) and \(e_{\mathcal Mz}\) are computed by replacing the four subcell values of \(\sigma_h\) and \(z_h\) on each macroelement by their macroelement averages and then taking the corresponding \(L^2\) errors against the exact fields.
	The observed errors \(e_u,e_\gamma,e_p\) are \(\mathcal O(h^2)\),
	whereas the stress and Darcy-velocity errors \(e_\sigma,e_z\)
	are \(\mathcal O(h)\). The macroelement-averaged stress and Darcy velocity exhibit second-order superconvergence.
	Since the same rates are retained for \(\lambda=10^6\), including the low-permeability case \(\kappa=10^{-6}\), the test indicates locking-free accuracy and no visible deterioration of the discretization error in these parameter regimes.

	\begin{figure}[htbp]
		\centering
		\setlength{\tabcolsep}{2pt}
		\begin{tabular}{ccc}
			\includegraphics[width=0.32\textwidth]{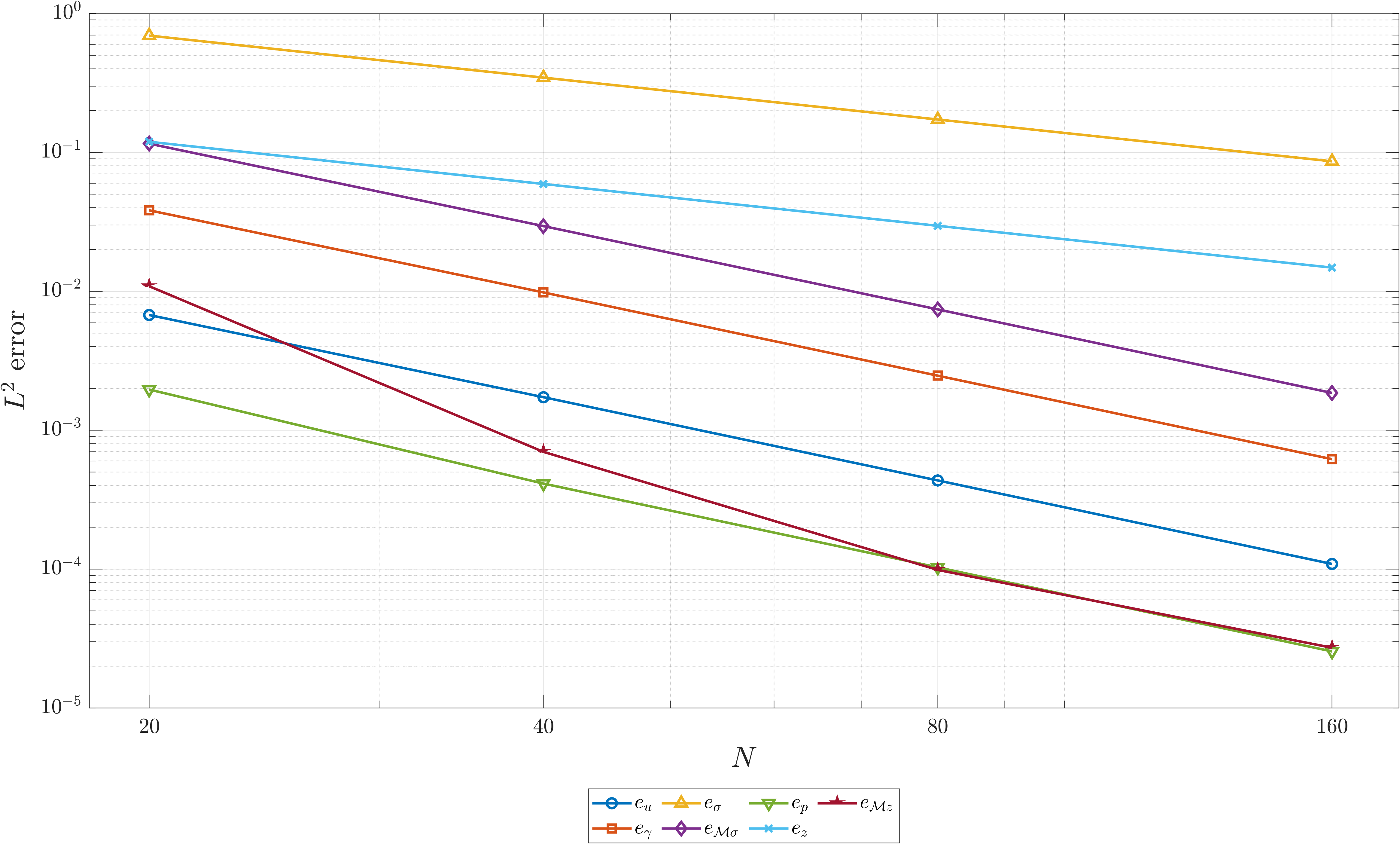}
			&
			\includegraphics[width=0.32\textwidth]{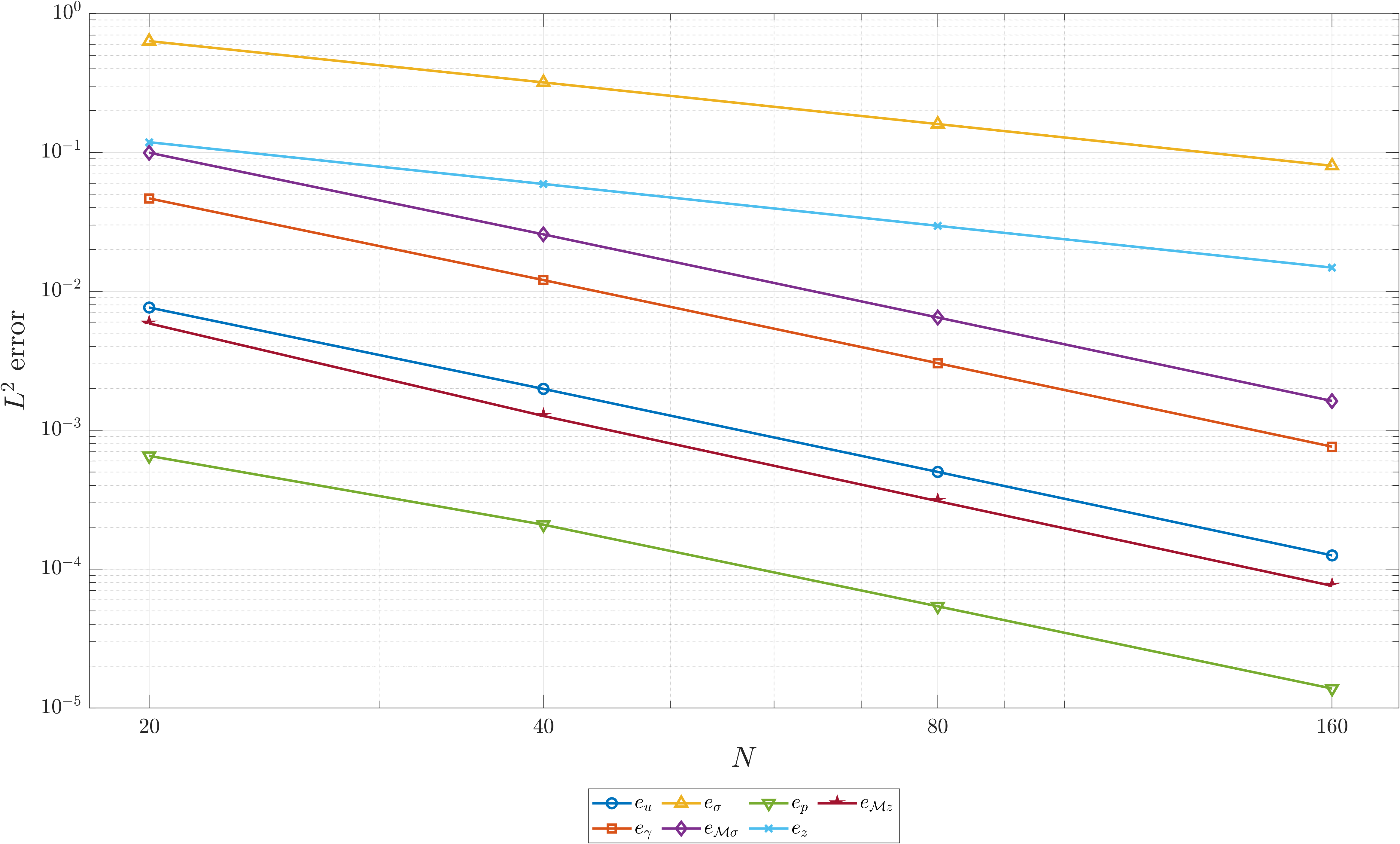}
			&
			\includegraphics[width=0.32\textwidth]{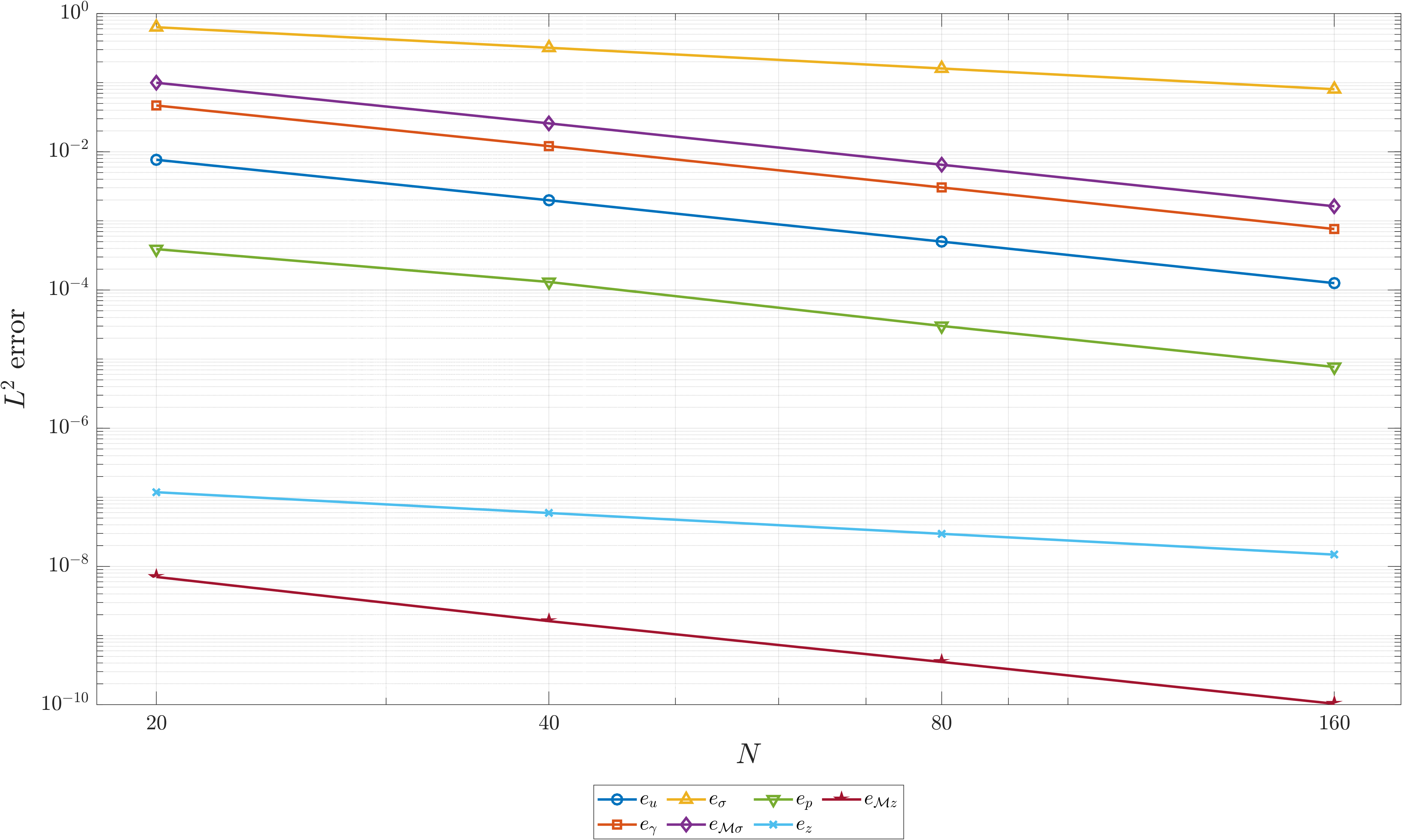}
		\end{tabular}
		\caption{\label{fig:example1-accuracy-convergence}Accuracy results for Example~\ref{example:accuracy-manufacture}. Left: \(\lambda=\kappa=1\); center: \(\lambda=10^6\), \(\kappa=1\); right: \(\lambda=10^6\), \(\kappa=10^{-6}\).}
	\end{figure}

	We next assess the parameter robustness of \(\mathcal P_{\rm bd}^{\Pmz}\) and \(\mathcal P_{\rm bd}^{\rm sc}\). Table~\ref{tab:example1-iterations} reports the GMRES iteration counts. Both preconditioners converge to the prescribed tolerance
	for all tested physical and discretization parameters. {The counts remain between 25 and 35 across the tested meshes and parameter combinations.} The Schur-consistent block \(\mathcal P_{\rm bd}^{\rm sc}\) never requires more outer iterations than the uncoupled block \(\mathcal P_{\rm bd}^{\Pmz}\), although the difference is generally small. To compare this iteration reduction with the actual computational cost, Table~\ref{tab:example1-one-step-time} reports the GMRES CPU time (in seconds) for one backward Euler step. The ratio \(t_{\rm sc}/t_{\Pmz}\) is the Schur-consistent solve time divided by the uncoupled solve time. The Schur-consistent preconditioner is about 4.7--10.4 times more expensive in the reported tests because each outer application contains an inner pressure solve. Consequently, the less expensive pressure block of \(\mathcal P_{\rm bd}^{\Pmz}\) gives the lower overall cost despite its slightly larger iteration counts.

	\begin{table}[htbp]
		\centering
		\setlength{\abovecaptionskip}{2pt}
		\setlength{\belowcaptionskip}{3pt}
		\caption{GMRES iterations for Example~\ref{example:accuracy-manufacture}, with \(c_0=1\).}
		\label{tab:example1-iterations}
		\papertablestyle
		\paramtablerules
		\begin{tabular}{@{}ccccc|ccc|ccc@{}}
			\toprule
			&& \multicolumn{3}{c}{\(\Delta t=10^{-2}\)} & \multicolumn{3}{c}{\(\Delta t=10^{-3}\)} & \multicolumn{3}{c}{\(\Delta t=10^{-4}\)} \\
			\cmidrule(lr){3-5}\cmidrule(lr){6-8}\cmidrule(l){9-11}
			\(\lambda\) & \(\kappa\) & \(N=40\) & \(80\) & \multicolumn{1}{c}{\(160\)} & \(N=40\) & \(80\) & \multicolumn{1}{c}{\(160\)} & \(N=40\) & \(80\) & \(160\) \\
			\midrule
			\(1\)    & \(1\)       & 33 (28) & 30 (27) & 28 (25) & 34 (29) & 33 (28) & 31 (27) & 33 (29) & 32 (28) & 32 (27) \\
			& \(10^{-2}\) & 33 (29) & 32 (28) & 32 (27) & 33 (29) & 32 (28) & 31 (27) & 32 (29) & 30 (26) & 30 (26) \\
			& \(10^{-6}\) & 32 (29) & 31 (26) & 29 (25) & 32 (29) & 31 (26) & 29 (25) & 32 (29) & 31 (26) & 29 (25) \\
			\hline
			\(10^3\) & \(1\)       & 32 (28) & 32 (28) & 29 (26) & 29 (29) & 31 (29) & 29 (27) & 29 (29) & 29 (29) & 27 (27) \\
			& \(10^{-2}\) & 29 (29) & 29 (29) & 27 (27) & 29 (29) & 29 (29) & 26 (26) & 29 (29) & 29 (28) & 26 (26) \\
			& \(10^{-6}\) & 29 (29) & 29 (28) & 26 (26) & 29 (29) & 29 (28) & 26 (26) & 29 (29) & 29 (28) & 26 (26) \\
			\hline
			\(10^6\) & \(1\)       & 34 (28) & 35 (28) & 35 (29) & 31 (28) & 33 (28) & 35 (29) & 28 (28) & 30 (28) & 32 (29) \\
			& \(10^{-2}\) & 28 (28) & 30 (28) & 32 (29) & 28 (28) & 28 (28) & 29 (29) & 28 (28) & 28 (28) & 29 (29) \\
			& \(10^{-6}\) & 28 (28) & 28 (28) & 29 (29) & 28 (28) & 28 (28) & 29 (29) & 28 (28) & 28 (28) & 29 (29) \\
			\bottomrule
		\end{tabular}
	\end{table}

	\begin{table}[htbp]
		\centering
		\caption{GMRES CPU time (s) and time ratio for Example~\ref{example:accuracy-manufacture}.}
		\label{tab:example1-one-step-time}
		\papertablestyle
		\begin{tabular}{@{}c*{4}{rrr}@{}}
			\toprule
			& \multicolumn{3}{c}{\(N=20\)} & \multicolumn{3}{c}{\(N=40\)} & \multicolumn{3}{c}{\(N=80\)} & \multicolumn{3}{c}{\(N=160\)} \\
			\cmidrule(lr){2-4}\cmidrule(lr){5-7}\cmidrule(lr){8-10}\cmidrule(l){11-13}
			Parameter	family & \({\mathcal P_{\rm bd}^{\Pmz}}\) & \({\mathcal P_{\rm bd}^{\rm sc}}\) & ratio
			& \({\mathcal P_{\rm bd}^{\Pmz}}\) & \({\mathcal P_{\rm bd}^{\rm sc}}\) & ratio
			& \({\mathcal P_{\rm bd}^{\Pmz}}\) & \({\mathcal P_{\rm bd}^{\rm sc}}\) & ratio
			& \({\mathcal P_{\rm bd}^{\Pmz}}\) & \({\mathcal P_{\rm bd}^{\rm sc}}\) & ratio \\
			\midrule
			1 & 0.070 & 0.566 & 8.14 & 0.286 & 2.084 & 7.30 & 1.219 & 12.732 & 10.44 & 5.742 & 52.492 & 9.14 \\
			2 & 0.123 & 0.675 & 5.49 & 0.611 & 2.885 & 4.72 & 2.835 & 16.279 & 5.74 & 15.006 & 95.905 & 6.39 \\
			\midrule
			\multicolumn{13}{@{}l}{Parameter family 1: \(\lambda=\kappa=c_0=1\), \(\Delta t=10^{-2}\);}
			\\[-1pt]
			\multicolumn{13}{@{}l}{Parameter family 2: \(\lambda=10^6\), \(\kappa=10^{-6}\), \(c_0=10^{-12}\), \(\Delta t=10^{-4}\).} \\
			\bottomrule
		\end{tabular}
	\end{table}

	We further test the critically balanced scaling \(c_0=\Kdr^{-1}\), with \(\mu=\alpha_B=1\), \(\kappa=10^{-6}\), and \(\Delta t=10^{-3}\). In 2D, \(\Kdr=\lambda+1\), so this choice gives \(\chi_{\lambda,c_0}=1\) while \(c_0\to0\) as \(\lambda\to\infty\). Table~\ref{tab:example1-critical-path-iterations} shows that all cases converge. The counts remain between 27 and 65, with no iteration blow-up as \(c_0\) approaches zero. These results support the robustness predicted for the balanced storage scaling.

	\begin{table}[htbp]
		\centering
		\caption{GMRES iterations along the path \(c_0=\Kdr^{-1}\) for Example~\ref{example:accuracy-manufacture}.}
		\label{tab:example1-critical-path-iterations}
		\papertablestyle
		\begin{tabular}{c*{7}{c}}
			\toprule
			\(N\backslash\lambda\)
			& \(10^0\) & \(10^1\) & \(10^2\) & \(10^3\)
			& \(10^4\) & \(10^5\) & \(10^6\) \\
			\midrule
			20  & 38 (30) & 37 (33) & 39 (37) & 42 (41) & 45 (43) & 46 (44) & 46 (44) \\
			40  & 36 (30) & 37 (33) & 42 (40) & 45 (43) & 49 (47) & 49 (47) & 49 (49) \\
			80  & 35 (29) & 34 (30) & 39 (36) & 44 (42) & 49 (47) & 55 (53) & 49 (49) \\
			160 & 33 (27) & 31 (29) & 35 (31) & 41 (40) & 52 (50) & 65 (62) & 53 (52) \\
			\bottomrule
		\end{tabular}
	\end{table}
\end{example}
\begin{example}[2D manufactured solution on perturbed quadrilateral meshes]\label{example:accuracy-manufacture-quad}
	\normalfont
	We repeat the manufactured-solution test of Example~\ref{example:accuracy-manufacture}, now on perturbed quadrilateral meshes. Let \((\xi_i,\eta_j)=(ih,jh)\), where \(h=1/N\), denote the vertices of the uniform logical grid. The physical vertices are defined by the smooth mapping
	\begin{equation*}x_{ij}=\xi_i+h^2\sin(2\pi\xi_i)\sin(2\pi\eta_j),
		\qquad
		y_{ij}=\eta_j+h^2\sin(2\pi\xi_i)\sin(2\pi\eta_j).
	\end{equation*}
	Thus, the boundary of \(\Omega=(0,1)^2\) is unchanged, whereas the interior vertices are displaced by \(\mathcal O(h^2)\), producing a perturbed family of general quadrilateral meshes. We use the same exact solution and fixed parameters \(\mu=\alpha_B=1\) as in Example~\ref{example:accuracy-manufacture}. For both practical preconditioners, the displacement block retains the general-quadrilateral two-grid solver.
	This experiment therefore isolates the effect of the \(\mathcal O(h^2)\) quadrilateral mesh perturbation on the performance of the stabilized formulation and its practical preconditioners. Table~\ref{tab:example2-quad-iterations-dt1e-2} shows that,
	across the tested meshes and parameter combinations, the uncoupled
	and Schur-consistent preconditioners require 26--104 and 25--53
	iterations, respectively. No systematic iteration growth is observed as the mesh is refined and the parameters vary.
	\begin{table}[htbp]
		\centering
		\setlength{\abovecaptionskip}{2pt}
		\setlength{\belowcaptionskip}{3pt}
		\caption{GMRES iterations for Example~\ref{example:accuracy-manufacture-quad}: \(\Delta t=10^{-2}\).}
		\label{tab:example2-quad-iterations-dt1e-2}
		\papertablestyle
		\paramtablerules
		\begin{tabular}{@{}ccccc|ccc|ccc@{}}
			\toprule
			&& \multicolumn{3}{c}{\(c_0=1\)} & \multicolumn{3}{c}{\(c_0=10^{-4}\)} & \multicolumn{3}{c}{\(c_0=10^{-8}\)} \\
			\cmidrule(lr){3-5}\cmidrule(lr){6-8}\cmidrule(l){9-11}
			\(\lambda\) & \(\kappa\) & \(N=40\) & \(80\) & \multicolumn{1}{c}{\(160\)} & \(N=40\) & \(80\) & \multicolumn{1}{c}{\(160\)} & \(N=40\) & \(80\) & \(160\) \\
			\midrule
			\(1\) & \(1\) & 33 (28) & 30 (27) & 28 (25) & 37 (29) & 34 (27) & 32 (26) & 37 (29) & 34 (27) & 32 (26) \\
			& \(10^{-2}\) & 33 (29) & 32 (28) & 32 (27) & 61 (39) & 60 (39) & 57 (36) & 61 (39) & 60 (39) & 57 (36) \\
			& \(10^{-6}\) & 33 (29) & 32 (27) & 30 (25) & 87 (36) & 85 (33) & 83 (31) & 104 (30) & 97 (33) & 93 (30) \\
			\hline
			\(10^3\) & \(1\) & 32 (28) & 32 (28) & 29 (26) & 32 (29) & 32 (29) & 29 (27) & 32 (29) & 32 (29) & 29 (27) \\
			& \(10^{-2}\) & 29 (29) & 29 (29) & 27 (27) & 40 (39) & 38 (38) & 37 (35) & 40 (39) & 38 (38) & 37 (35) \\
			& \(10^{-6}\) & 29 (29) & 29 (28) & 26 (26) & 55 (49) & 55 (48) & 52 (46) & 66 (49) & 65 (48) & 61 (45) \\
			\hline
			\(10^6\) & \(1\) & 37 (30) & 35 (28) & 35 (29) & 38 (30) & 36 (28) & 36 (29) & 38 (30) & 36 (28) & 36 (29) \\
			& \(10^{-2}\) & 30 (30) & 30 (28) & 32 (29) & 38 (32) & 36 (30) & 36 (31) & 38 (32) & 36 (30) & 36 (31) \\
			& \(10^{-6}\) & 30 (30) & 28 (28) & 29 (29) & 37 (37) & 35 (35) & 36 (36) & 57 (53) & 53 (49) & 52 (49) \\
			\bottomrule
		\end{tabular}
	\end{table}

\end{example}

\begin{example}[Barry--Mercer localized-source problem]
	\label{example:barry-mercer-like}
	\normalfont
	We then consider a 2D Barry--Mercer benchmark to illustrate the pressure behavior in a low-permeability, small-time-step regime. The domain is \(\Omega=(0,1)^2\), the initial data are zero, and the point source at \((x_0,y_0)=(1/4,1/4)\) is represented by \(q_h^n=Q(t_n)|M_0|^{-1}\chi_{M_0}\), where \(Q(t)=\cos(2\pi t)\) and \(M_0\) is the source cell.  Homogeneous boundary data are imposed for displacement and pressure. We let \(E\) and \(\nu\) denote Young's modulus and Poisson's ratio, respectively. Under the plane-strain assumption, the Lam\'e coefficients are \(\lambda=E\nu/((1+\nu)(1-2\nu))\) and \(\mu=E/(2(1+\nu))\). Following the challenging parameter regime used in \cite{PeDeLaRivaEtAl2025}, we take \(E=10^5\), \(\nu=0.1\), \(\alpha_B=1\), \(c_0=10^{-8}\), \(\Delta t=10^{-4}\), and \(\kappa=10^{-6}\).
	Figure~\ref{fig:barry-mercer-oscillation} shows the computed pressure at \(T=10^{-4}\) on the \(64\times64\) mesh. The left panel gives the full pressure surface, whereas the right panel is the horizontal cut through the source row. Although the low permeability produces a sharply localized pressure peak, the pressure cut exhibits a clean single-peak profile followed by a smooth decay. The computed pressure remains positive,
	and no visible cell-to-cell saw-tooth oscillations occur near or away from the source.
	This indicates that the proposed control-volume formulation can effectively suppress nonphysical pressure oscillations without a separate pressure stabilization term in the flow equation.

	\begin{figure}[htbp]
		\centering
		\setlength{\tabcolsep}{3pt}
		\begin{tabular}{cc}
			\includegraphics[width=0.47\textwidth]{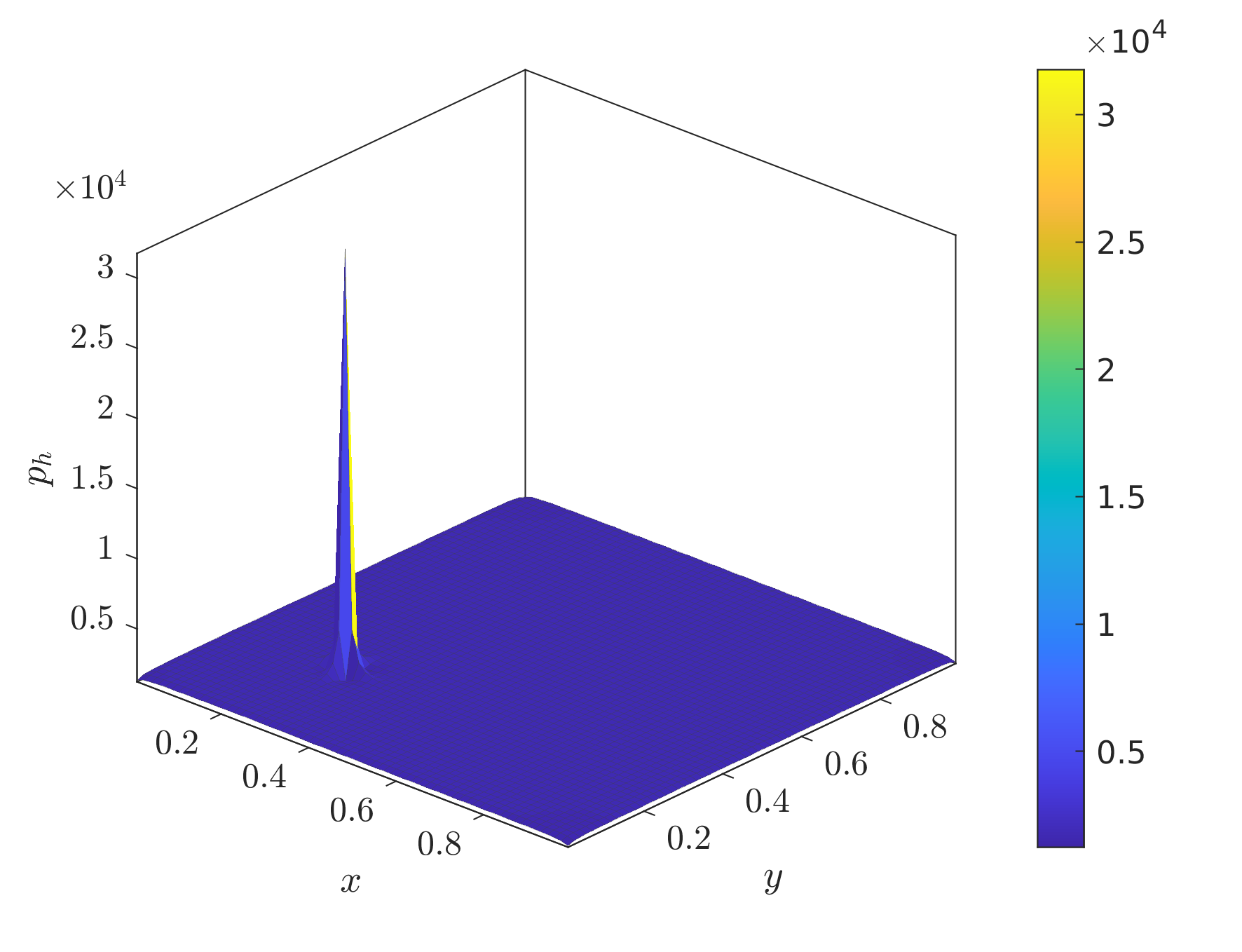}
			&
			\includegraphics[width=0.47\textwidth]{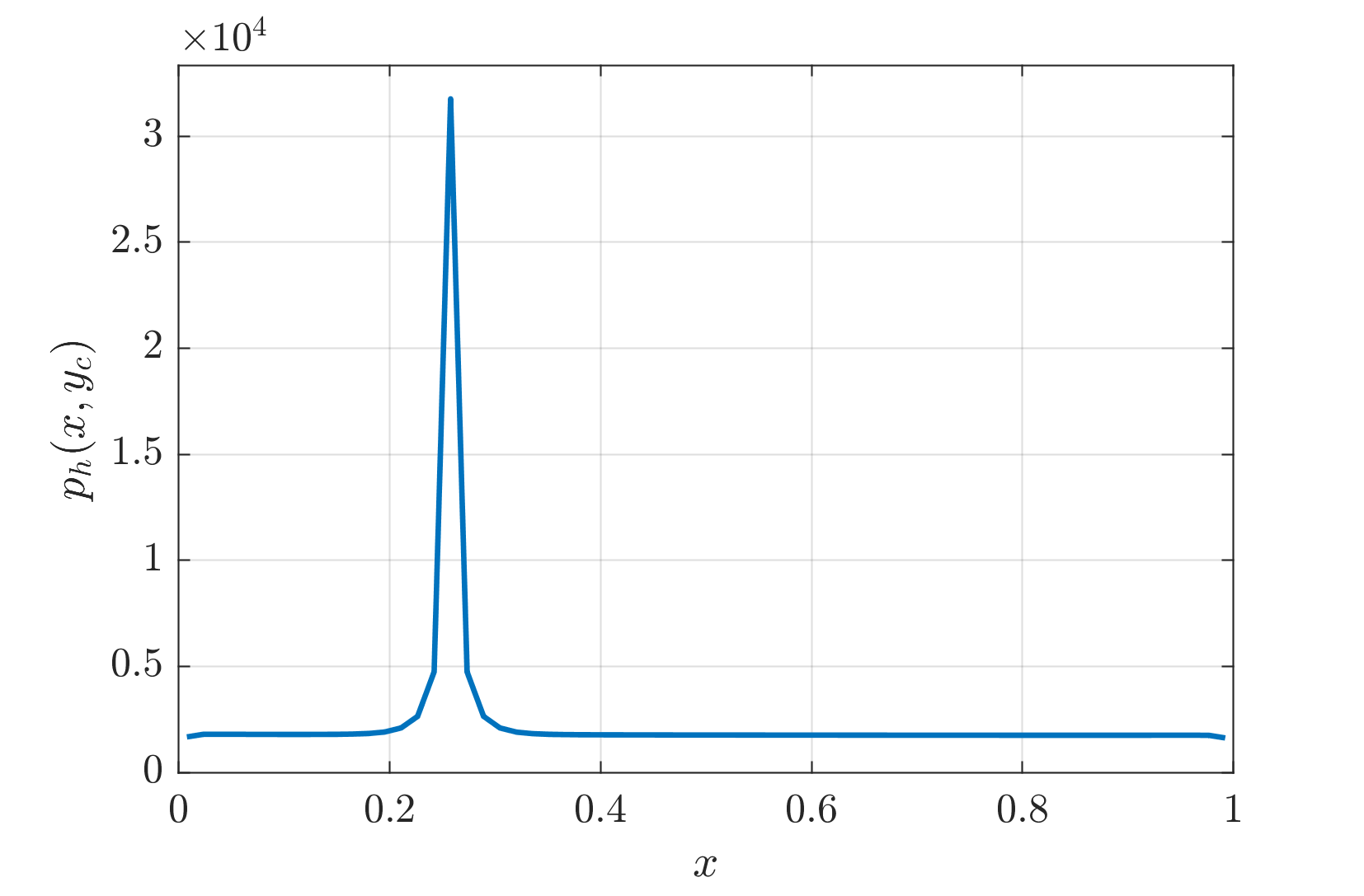}
			\\[-2pt]
			(a) Pressure surface & (b) Pressure cut through the source row
		\end{tabular}
		\caption{Pressure behavior for Example~\ref{example:barry-mercer-like} on the \(64\times64\) mesh.}
		\label{fig:barry-mercer-oscillation}
	\end{figure}

	We also use this localized-source problem to probe the parameter behavior of the two block-diagonal preconditioners. We keep \(E=10^5\), \(\alpha_B=1\), and \(\Delta t=T=10^{-4}\), and choose \(\nu\in\{0.1,0.4,0.49\}\).
	For the comparison, we fix \(\kappa=10^{-6}\) and \(c_0=10^{-8}\). Table~\ref{tab:barry-mercer-parameter-iterations} reports the reduced-system GMRES counts. Both preconditioners converge to the prescribed tolerance.
	As \(\lambda\) increases over the tested range, the uncoupled
	iteration counts decrease from 96 to 46 at \(N=40\)
	and from 84 to 40 at \(N=160\).
	The Schur-consistent block-diagonal preconditioner remains particularly stable, requiring only 30--35 outer iterations throughout the reported parameter range.

	\begin{table}[htbp]
		\centering
		\caption{GMRES iterations for Example~\ref{example:barry-mercer-like}.}
		\label{tab:barry-mercer-parameter-iterations}
		\papertablestyle
		\begin{tabular}{@{}cccc@{}}
			\toprule
			\(\lambda\) & \(N=40\) & \(N=80\) & \(N=160\) \\
			\midrule
			\(1.1364\times10^4\) & 96 (30) & 91 (33) & 84 (35) \\
			\(1.4286\times10^5\) & 65 (31) & 58 (32) & 56 (32) \\
			\(1.6443\times10^6\) & 46 (34) & 42 (33) & 40 (31) \\
			\bottomrule
		\end{tabular}
	\end{table}
\end{example}

\begin{example}[3D manufactured solution on structured cubical meshes]\label{example:3d-manufacture}
	\normalfont
	On \(\Omega=(0,1)^3\), we take \(p=e^{-t}\sin(\pi x)\sin(\pi y)\sin(\pi z)\),
\(u_1=p\), \(u_2=2p\cos(\pi x)\), and \(u_3=2p\cos(\pi y)\).
The initial pressure and stress are obtained from these fields at \(t=0\).
We set \(\mu=\alpha_B=1\) and \(\Delta t=10^{-3}\).
	Table~\ref{tab:example3d-parameter-iterations} reports iteration counts for the block preconditioners applied to the three-field system as the physical parameters $\lambda$, $\kappa$, and $c_0$ vary.
	Both preconditioners reach the prescribed tolerance. Thus, neither large \(\lambda\) nor small \(\kappa\) causes an iteration blow-up. Even for \(c_0=0\), \(\mathcal P_{\rm bd}^{\Pmz}\) still converges in 15--57 iterations.

	\begin{table}[htbp]
		\centering
		\setlength{\abovecaptionskip}{2pt}
		\setlength{\belowcaptionskip}{3pt}
		\caption{GMRES iterations for Example~\ref{example:3d-manufacture}.}
		\label{tab:example3d-parameter-iterations}
		\papertablestyle
		\paramtablerules
		\begin{tabular}{@{}cccccc|cccc|cccc@{}}
			\toprule
			&& \multicolumn{4}{c}{\(c_0=1\)} & \multicolumn{4}{c}{\(c_0=10^{-4}\)} & \multicolumn{4}{c}{\(c_0=0\)} \\
			\cmidrule(lr){3-6}\cmidrule(lr){7-10}\cmidrule(l){11-14}
			\(\lambda\) & \(\kappa\) & \(N=16\) & 32 & 48 & \multicolumn{1}{c}{64} & \(N=16\) & 32 & 48 & \multicolumn{1}{c}{64} & \(N=16\) & 32 & 48 & 64 \\
			\midrule
			\(10^2\) & \(10^{-2}\) & 20 (20) & 20 (20) & 19 (19) & 19 (19) & 31 (27) & 31 (27) & 32 (26) & 31 (26) & 31 (27) & 31 (27) & 32 (26) & 31 (26) \\
			& \(10^{-4}\) & 20 (20) & 20 (20) & 19 (19) & 19 (19) & 40 (27) & 38 (27) & 38 (26) & 36 (26) & 44 (27) & 43 (27) & 41 (26) & 39 (26) \\
			& \(10^{-6}\) & 20 (20) & 20 (20) & 19 (19) & 19 (19) & 41 (27) & 41 (27) & 38 (25) & 38 (24) & 57 (31) & 57 (31) & 52 (29) & 26 (28) \\
			\hline
			\(10^4\) & \(10^{-2}\) & 20 (20) & 20 (20) & 22 (22) & 23 (23) & 22 (16) & 19 (16) & 19 (14) & 18 (14) & 24 (16) & 21 (15) & 20 (14) & 20 (14) \\
			& \(10^{-4}\) & 20 (20) & 20 (20) & 22 (22) & 23 (23) & 21 (21) & 20 (18) & 18 (16) & 17 (15) & 28 (21) & 26 (20) & 24 (18) & 23 (17) \\
			& \(10^{-6}\) & 20 (20) & 20 (20) & 22 (22) & 23 (23) & 21 (21) & 18 (16) & 17 (15) & 15 (13) & 21 (21) & 18 (18) & 17 (17) & 15 (15) \\
			\hline
			\(10^6\) & \(10^{-2}\) & 19 (19) & 22 (22) & 24 (24) & 26 (24) & 28 (18) & 27 (16) & 27 (16) & 26 (16) & 29 (18) & 27 (16) & 27 (16) & 26 (16) \\
			& \(10^{-4}\) & 19 (19) & 22 (22) & 24 (24) & 24 (24) & 17 (17) & 18 (15) & 18 (15) & 18 (15) & 23 (16) & 20 (14) & 18 (14) & 17 (13) \\
			& \(10^{-6}\) & 19 (19) & 22 (22) & 24 (24) & 24 (24) & 17 (17) & 14 (14) & 14 (14) & 13 (13) & 21 (21) & 20 (18) & 20 (15) & 19 (15) \\
			\bottomrule
		\end{tabular}
	\end{table}

	Table~\ref{tab:example3d-one-step-time} reports the GMRES CPU times for two selected parameter families. The outer iteration counts of \(\mathcal P_{\rm bd}^{\rm sc}\) do not exceed those of \(\mathcal P_{\rm bd}^{\Pmz}\) in these tests. Nevertheless, \(\mathcal P_{\rm bd}^{\Pmz}\) requires less CPU time on every tested mesh, with \(t_{\rm sc}/t_{\Pmz}\) ranging from 5.12 to 9.54.

	\begin{table}[htbp]
		\centering
		\caption{GMRES CPU time (s) and time ratio for Example~\ref{example:3d-manufacture}.}
		\label{tab:example3d-one-step-time}
		\papertablestyle
		\begin{tabular}{@{}c*{4}{rrr}@{}}
			\toprule
			& \multicolumn{3}{c}{\(N=16\)} & \multicolumn{3}{c}{\(N=32\)} & \multicolumn{3}{c}{\(N=48\)} & \multicolumn{3}{c}{\(N=64\)} \\
			\cmidrule(lr){2-4}\cmidrule(lr){5-7}\cmidrule(lr){8-10}\cmidrule(l){11-13}
			Parameter family & \({\mathcal P_{\rm bd}^{\Pmz}}\) & \({\mathcal P_{\rm bd}^{\rm sc}}\) & ratio
			& \({\mathcal P_{\rm bd}^{\Pmz}}\) & \({\mathcal P_{\rm bd}^{\rm sc}}\) & ratio
			& \({\mathcal P_{\rm bd}^{\Pmz}}\) & \({\mathcal P_{\rm bd}^{\rm sc}}\) & ratio
			& \({\mathcal P_{\rm bd}^{\Pmz}}\) & \({\mathcal P_{\rm bd}^{\rm sc}}\) & ratio \\
			\midrule
			1 & 7.0 & 35.6 & 5.12 & 72.2 & 398.9 & 5.52 & 338.1 & 2025.5 & 5.99 & 1258.6 & 8697.4 & 6.91 \\
			2 & 7.2 & 42.3 & 5.89 & 72.4 & 494.6 & 6.83 & 319.2 & 2193.1 & 6.87 & 1129.7 & 10781.0 & 9.54 \\
			\midrule
			\multicolumn{13}{@{}l}{Parameter family 1: \(\lambda=10^2\), \(\kappa=10^{-2}\), \(c_0=1\), \(\Delta t=10^{-3}\);}\\[-1pt]
			\multicolumn{13}{@{}l}{Parameter family 2: \(\lambda=10^6\), \(\kappa=10^{-6}\), \(c_0=0\), \(\Delta t=10^{-3}\).} \\
			\bottomrule
		\end{tabular}
	\end{table}

\end{example}

\begin{example}[2D and 3D footing problems]\label{example:footing}
	\normalfont
	The footing problem is a classical Biot benchmark \cite{GasparLisbonaOosterlee2008}. We follow the 2D initial--boundary value setting in \cite{HenriquezEtAl2025Biot} on \(\Omega=(-50,50)\times(0,75)\), with downward traction \(\underline\sigma\bs n=(0,-\sigma_0)^{\top}\), \(\sigma_0=10^4\), on \(\Gamma_1=\{(x,75):|x|\leq50/3\}\). The remaining top boundary is traction-free, the sides and bottom are fixed, and \(p=0\) on \(\partial\Omega\). The body force, fluid source, and initial data are zero.
	For the representative solution in Figure~\ref{fig:footing2d-physical-fields}, we use the benchmark parameters \(E=3\times10^4\), \(\nu=0.4995\), \(\kappa=10^{-4}\), \(c_0=10^{-3}\), and \(\alpha_B=0.1\), together with \(N=60\), \(\Delta t=1\), and final time \(T=50\). The snapshots are shown at \(t=1,10,25,50\). The black segment on the top boundary marks the loaded region \(\Gamma_1\). The pressure remains smooth and symmetric, with its largest values beneath \(\Gamma_1\); as time advances, its peak decreases and the low-pressure layers adjacent to the drained boundary extend inward. The logarithmic flux plots reveal drainage over several orders of magnitude, with the strongest early flow occurring near the boundary and the two endpoints of \(\Gamma_1\). These high-flux regions steadily weaken, with \(\max|\bs z_h|\) decreasing from \(1.34\times10^{-5}\) at \(t=1\) to \(2.19\times10^{-6}\) at \(t=50\), and no visible oscillatory pattern develops.
	\begin{figure}[htbp]
		\centering
		\includegraphics[width=\textwidth]{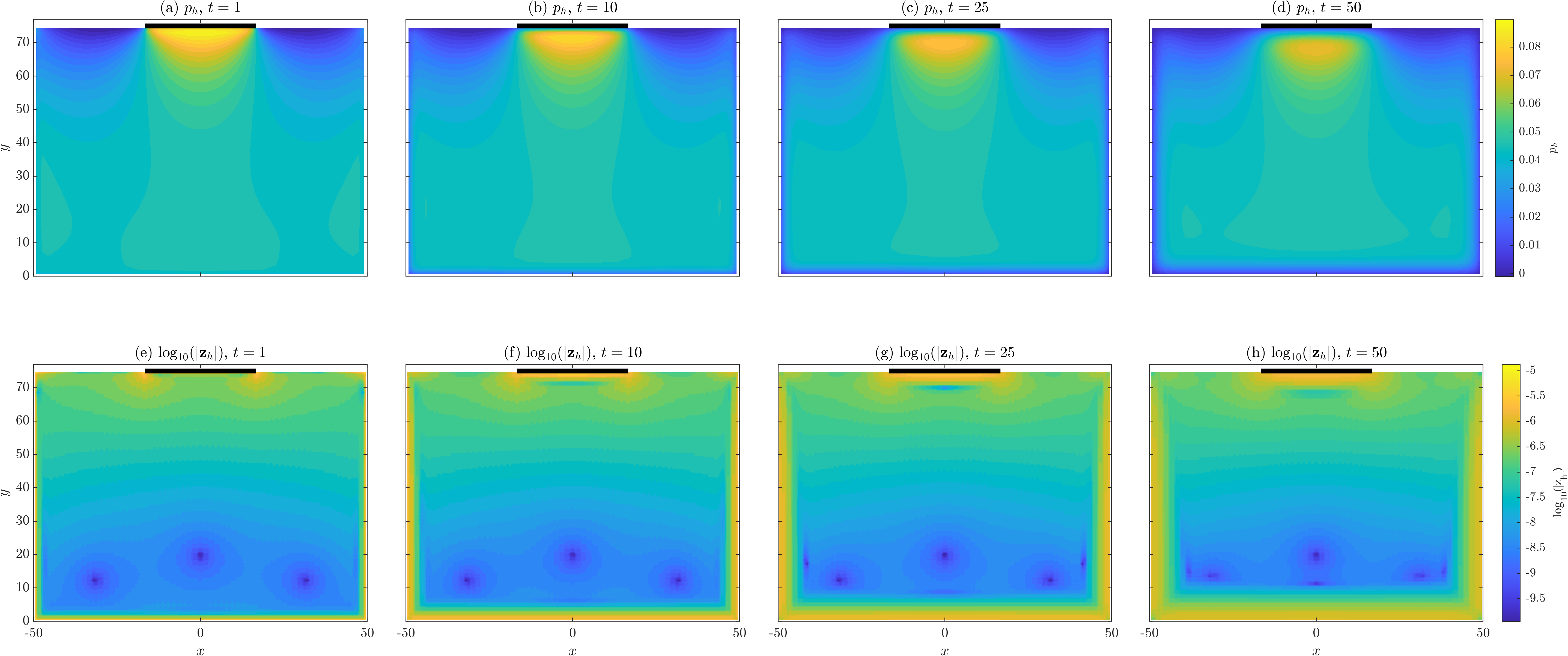}
		\caption{Evolution of the pressure \(p_h\) (top row) and logarithmic Darcy-velocity magnitude \(\log_{10}|\bs z_h|\) (bottom row) for the 2D setting of Example~\ref{example:footing} at \(t=1,10,25,50\). }
		\label{fig:footing2d-physical-fields}
	\end{figure}

	To test parameter robustness, we take \(\Delta t=10^{-2}\) and fix \(\mu=1\) and \(\alpha_B=0.1\). Table~\ref{tab:footing2d-parameter-iterations} reports the GMRES iteration counts.
	{The} two preconditioners remain effective as \(\lambda\), \(\kappa\), and \(c_0\) vary over the reported ranges, including the limiting case \(c_0=0\); no parameter choice causes an iteration blow-up.
	\begin{table}[htbp]
		\centering
		\setlength{\abovecaptionskip}{2pt}
		\setlength{\belowcaptionskip}{3pt}
		\caption{GMRES iterations for the 2D setting of Example~\ref{example:footing}.}
		\label{tab:footing2d-parameter-iterations}
		\papertablestyle
		\paramtablerules
		\begin{tabular}{@{}ccccc|ccc|ccc@{}}
			\toprule
			&& \multicolumn{3}{c}{\(c_0=1\)} & \multicolumn{3}{c}{\(c_0=10^{-4}\)} & \multicolumn{3}{c}{\(c_0=0\)} \\
			\cmidrule(lr){3-5}\cmidrule(lr){6-8}\cmidrule(l){9-11}
			\(\lambda\) & \(\kappa\) & \(N=40\) & 80 & \multicolumn{1}{c}{160} & \(N=40\) & 80 & \multicolumn{1}{c}{160} & \(N=40\) & 80 & 160 \\
			\midrule
			\(1\) & \(10^{-2}\) & 34 (33) & 35 (33) & 35 (33) & 88 (58) & 92 (60) & 94 (61) & 90 (59) & 94 (60) & 96 (61) \\
			& \(10^{-4}\) & 34 (33) & 35 (33) & 35 (33) & 88 (58) & 94 (60) & 97 (59) & 90 (59) & 96 (60) & 99 (61) \\
			& \(10^{-6}\) & 34 (33) & 35 (33) & 35 (33) & 88 (58) & 94 (60) & 97 (59) & 90 (59) & 96 (61) & 99 (61) \\
			\hline
			\(10^4\) & \(10^{-2}\) & 27 (27) & 28 (28) & 29 (29) & 33 (33) & 35 (35) & 35 (35) & 50 (50) & 54 (54) & 54 (54) \\
			& \(10^{-4}\) & 27 (27) & 28 (28) & 29 (29) & 33 (33) & 35 (35) & 35 (35) & 60 (60) & 66 (66) & 66 (66) \\
			& \(10^{-6}\) & 27 (27) & 28 (28) & 29 (29) & 33 (33) & 35 (35) & 35 (35) & 61 (61) & 66 (66) & 67 (67) \\
			\hline
			\(10^6\) & \(10^{-2}\) & 27 (27) & 29 (29) & 41 (41) & 28 (28) & 32 (32) & 46 (46) & 33 (33) & 37 (37) & 54 (52) \\
			& \(10^{-4}\) & 27 (27) & 29 (29) & 41 (41) & 28 (28) & 32 (32) & 46 (46) & 50 (50) & 54 (54) & 75 (77) \\
			& \(10^{-6}\) & 27 (27) & 29 (29) & 41 (42) & 28 (28) & 32 (32) & 46 (46) & 61 (60) & 67 (66) & 99 (95) \\
			\bottomrule
		\end{tabular}
	\end{table}

	We then use the 3D footing geometry of \cite{HenriquezEtAl2025Biot}, with \(\Omega=(-32,32)^2\times(0,64)\) and downward traction \(\underline\sigma\bs n=(0,0,-\sigma_0)^{\top}\) on \(\Gamma_1=\{(x,y,64):|x|\leq16,\ |y|\leq16\}\). The rest of the top face is traction-free, and the sides and bottom are fixed. The drainage conditions, initial data, body force, and fluid source are the same as in 2D.
	We use uniform structured cubical meshes with \(N^3\) cells and fix \(c_0=0.5\), \(\mu=1\), \(\alpha_B=0.5\), and \(\sigma_0=0.1\). Table~\ref{tab:footing3d-parameter-iterations} reports the GMRES iteration counts. The iteration counts range from 15 to 37 for \(\mathcal P_{\rm bd}^{\Pmz}\) and from 15 to 32 for \(\mathcal P_{\rm bd}^{\rm sc}\), with no deterioration as \(\kappa\) and \(\Delta t\) decrease over the tested ranges. Thus, both preconditioners remain effective under mesh refinement and the reported parameter variations.
	\begin{table}[htbp]
		\centering
		\setlength{\abovecaptionskip}{2pt}
		\setlength{\belowcaptionskip}{3pt}
		\caption{GMRES iterations for the 3D setting of Example~\ref{example:footing}.}
		\label{tab:footing3d-parameter-iterations}
		\papertablestyle
		\paramtablerules
		\begin{tabular}{@{}ccc|cccccc@{}}
			\toprule
			\(\lambda\) & \(\kappa\) & \(\Delta t\) & \(N=8\) & 16 & 24 & 32 & 48 & 64 \\
			\midrule
			\(1\) & \(10^{-4}\) & \(10^{-2}\) & 28 (23) & 33 (28) & 36 (31) & 36 (31) & 37 (31) & 37 (32) \\
			&  & \(10^{-3}\) & 28 (23) & 33 (28) & 36 (31) & 36 (31) & 37 (31) & 37 (32) \\
			& \(10^{-6}\) & \(10^{-2}\) & 28 (23) & 33 (28) & 36 (31) & 36 (31) & 37 (31) & 37 (32) \\
			&  & \(10^{-3}\) & 28 (23) & 33 (28) & 36 (31) & 36 (31) & 37 (31) & 37 (32) \\
			\hline
			\(10^4\) & \(10^{-4}\) & \(10^{-2}\) & 17 (17) & 20 (20) & 20 (20) & 21 (21) & 21 (21) & 23 (23) \\
			&  & \(10^{-3}\) & 17 (17) & 20 (20) & 20 (20) & 21 (21) & 21 (21) & 23 (23) \\
			& \(10^{-6}\) & \(10^{-2}\) & 17 (17) & 20 (20) & 20 (20) & 21 (21) & 21 (21) & 23 (23) \\
			&  & \(10^{-3}\) & 17 (17) & 20 (20) & 20 (20) & 21 (21) & 21 (21) & 23 (23) \\
			\hline
			\(10^6\) & \(10^{-4}\) & \(10^{-2}\) & 15 (15) & 18 (18) & 19 (19) & 19 (19) & 22 (22) & 24 (24) \\
			&  & \(10^{-3}\) & 15 (15) & 18 (18) & 19 (19) & 19 (19) & 22 (22) & 24 (24) \\
			& \(10^{-6}\) & \(10^{-2}\) & 15 (15) & 18 (18) & 19 (19) & 19 (19) & 22 (22) & 24 (24) \\
			&  & \(10^{-3}\) & 15 (15) & 18 (18) & 19 (19) & 19 (19) & 22 (22) & 24 (24) \\
			\bottomrule
		\end{tabular}
	\end{table}

\end{example}

\section{Conclusions}
\label{sec:conclusion}

We have developed and analyzed block preconditioners for a control-volume discretization of the quasi-static Biot equations. Starting from a five-field formulation, the stress and scaled Darcy velocity are eliminated by independent vertex-local solves, resulting in a symmetric condensed three-field system. We established the well-posedness of the fully discrete scheme and proved that the coupled mechanics block is uniformly spectrally equivalent to its displacement--rotation block-diagonal approximation.

Two classes of pressure preconditioners were then studied. For the inexpensive uncoupled strategy, we identified an exact discrete pressure-control condition and derived a sufficient condition. In particular, the analysis permits the critical scaling \(c_0\asymp \Kdr^{-1}\), so that the storage coefficient may tend to zero as the material approaches the nearly incompressible limit. For the Schur-consistent strategy, we introduced a computable pressure block and proved its uniform equivalence to the exact pressure Schur complement. This gives parameter-independent spectral bounds for Schur-consistent block-diagonal preconditioning and uniform weighted field-of-values estimates for Schur-consistent lower-triangular preconditioning.
Practical implementations of the block preconditioners were constructed using a vertex-patch two-grid solver for displacement, a macrocell mass inverse for the stabilized rotation block, and scalar two-grid solvers for pressure. The numerical experiments in 2D and 3D verify the expected accuracy of the control-volume discretization and demonstrate stable performance under mesh refinement and broad variations of the physical parameters.

\bibliographystyle{siamplain-nourl}
\bibliography{pre2026biotref}

\end{document}